\RequirePackage{ifpdf}
\ifpdf 
\documentclass[pdftex]{sigma}
\else
\documentclass{sigma}
\fi

\usepackage[all]{xy}

\numberwithin{equation}{section}

\newtheorem{Theorem}{Theorem}[section]
\newtheorem{Proposition}[Theorem]{Proposition}
\newtheorem{Lemma}[Theorem]{Lemma}
\newtheorem{Corollary}[Theorem]{Corollary}

\theoremstyle{definition}
\newtheorem{Definition}[Theorem]{Definition}
\newtheorem{Example}[Theorem]{Example}
\newtheorem{Remark}[Theorem]{Remark}

\begin{document}

\allowdisplaybreaks

\renewcommand{\PaperNumber}{102}

\FirstPageHeading

\ShortArticleName{Classical and Quantum Dilogarithm Identities}

\ArticleName{Classical and Quantum Dilogarithm Identities}

\Author{Rinat M.~KASHAEV~$^\dag$ and Tomoki NAKANISHI~$^\ddag$}

\AuthorNameForHeading{R.M.~Kashaev and T.~Nakanishi}

\Address{$^\dag$~Section de Math\'ematiques, Universit\'e de Gen\`eve,\\
\hphantom{$^\dag$}~2-4 rue du Li\`evre, Case postale 64, 1211 Gen\`eve 4, Switzerland}
\EmailD{\href{mailto:Rinat.Kashaev@unige.ch}{Rinat.Kashaev@unige.ch}}

\Address{$^\ddag$~Graduate School of Mathematics, Nagoya University,
Nagoya, 464-8604, Japan}
\EmailD{\href{mailto:nakanisi@math.nagoya-u.ac.jp}{nakanisi@math.nagoya-u.ac.jp}}

\ArticleDates{Received May 03, 2011, in f\/inal form October 26, 2011;  Published online November 01, 2011}

\Abstract{Using the quantum cluster algebra formalism of Fock and Goncharov,
we present several forms of quantum dilogarithm identities associated with periodicities in quantum cluster algebras,
namely, the tropical, universal, and local forms.
We then demonstrate how classical dilogarithm identities naturally
emerge from quantum dilogarithm identities in local form
 in the semiclassical limit by applying the saddle point method.}


\Classification{13F60}

\section{Introduction}

\subsection{Pentagon relations}

The {\em Euler dilogarithm\/}  $\mathrm{Li_2}(x)$
and its variant,
the {\em Rogers dilogarithm\/} $L(x)$
have appeared in several branches of
mathematics (e.g.,
 \cite{Lewin81,Kirillov95,Zagier07}).
See \eqref{eq:L00} and \eqref{eq:L0} for the def\/inition.
The most important property of the functions is
the {\em pentagon relation\/}. For $L(x)$, it takes the following form
\begin{gather}
\label{eq:Lpent}
L(x)+L(y)=
L\left(\frac{x(1-y)}{1-xy}\right)
+
L(xy)+
L\left(\frac{y(1-x)}{1-xy}\right),
\qquad 0\leq x,y\leq 1.
\end{gather}

The {\em quantum dilogarithm\/}  appears also
in several branches of
mathematics, e.g.,
 discrete quantum systems
\cite{Bazhanov93,Faddeev93,Faddeev94,Faddeev94b,Faddeev95,Faddeev01,
Bazhanov07,Bazhanov08,Kashaev08},
 hyperbolic geometry and Teichm\"uller theory
\cite{Kashaev97,Kashaev98,Fock07,Goncharov07},
quantum topology
\cite{Kashaev94,Baseilhac04},
Donaldson--Thomas
invariants
\cite{Kontsevich08,Kontsevich09,Kontsevich10,Keller11,Nagao11,Nagao11b},
string theory
\cite{Gaiotto08,Gaiotto10,Cecotti10},
 representation theory of algebras \cite{Reineke08}, etc.,
and it accumulates much attention recently.

Actually, there are at least two variants of the quantum dilogarithm.

The f\/irst one $\mathbf{\Psi}_q(x)$, where $q$ is a parameter,
 is simply called
 the {\em quantum dilogarithm\/} here.
See \eqref{eq:psi1} for the def\/inition.
The study of the function
as `quantum exponential' goes back to
\cite{Schutzenberger53},
but the recognition
 as  `quantum dilogarithm' was made more recently \cite{Faddeev93,Faddeev94}.
The following properties explain why
 it is considered as a quantum analogue of  the dilogarithm  \cite{Faddeev93,Faddeev94,Kashaev04}.

{\em
$(a)$ Asymptotic behavior:
  In  the limit $q \rightarrow 1^{-}$,
\begin{gather}
\label{eq:asym3}
\mathbf{\Psi}_q(x) \sim
\exp \left( -\frac{\mathrm{Li}_2(-x)}{2 \log q }
\right).
\end{gather}

$(b)$ Pentagon relation: If $UV=q^2VU$, then{\samepage
\begin{gather}
\label{eq:pent2}
\mathbf{\Psi}_q(U)\mathbf{\Psi}_q(V)
=
\mathbf{\Psi}_q(V)\mathbf{\Psi}_q\big(q^{-1}UV\big)\mathbf{\Psi}_q(U).
\end{gather}
Moreover, in the limit $q\rightarrow 1^{-}$, the relation
\eqref{eq:pent2} reduces to the   relation
\eqref{eq:Lpent}.}
}

The second variant of the quantum dilogarithm $\mathbf{\Phi}_b(z)$,
where $b$ is a parameter,
was introduced by \cite{Faddeev94b,Faddeev95}.
Here we call it {\em Faddeev's quantum dilogarithm}
(also known as the {\em noncompact quantum
dilogarithm}). See \eqref{eq:phi0}
for the def\/inition.
The function $\mathbf{\Phi}_b(z)$ also satisf\/ies properties
parallel to the ones for $\mathbf{\Psi}_q(x)$
\cite{Faddeev94b,Faddeev95,Woronowicz00, Faddeev01}.

{\em
$(a)$ Asymptotic behavior:  In  the limit $b \rightarrow 0$,
\begin{gather}
\label{eq:asym4}
\mathbf{\Phi}_b\left(\frac{z}{2\pi b}\right) \sim
\exp \left( -\frac{\mathrm{Li}_2(-e^{ z})}{2 \pi b^2 \sqrt{-1} }
\right).
\end{gather}

$(b)$ Pentagon relation: If  $[\hat P,\hat Q]
=(2\pi \sqrt{-1})^{-1}$, then
\begin{gather}
\label{eq:pent4}
\mathbf{\Phi}_b(\hat Q)\mathbf{\Phi}_b(\hat P)
=
\mathbf{\Phi}_b(\hat P)
\mathbf{\Phi}_b(\hat P+ \hat Q)
\mathbf{\Phi}_b(\hat Q).
\end{gather}
Moreover, in the limit $b\rightarrow 0$, the relation
\eqref{eq:pent4} reduces to the   relation
\eqref{eq:Lpent}.
}

Despite the appearance of the Euler dilogarithm
$\mathrm{Li}_2(x)$ in
\eqref{eq:asym3} and \eqref{eq:asym4},
  we have the Rogers dilogarithm $L(x)$
in \eqref{eq:Lpent} when we take the limits
of \eqref{eq:pent2}
and \eqref{eq:pent4}.
Namely, the limits of \eqref{eq:pent2}
and \eqref{eq:pent4} are not so trivial as termwise limit.
The two functions $L(x)$ and $\mathrm{Li}_2(x)$ dif\/fer by
logarithms (see \eqref{eq:LL1} and \eqref{eq:LL2}),
and the noncommutativity of $U$, $V$ and $P$, $Q$
`magically' turns $\mathrm{Li}_2(x)$ into~$L(x)$.
To clarify this phenomenon in a (much) wider situation
is the main theme of this paper.

\subsection{Classical and quantum dilogarithm identities from cluster algebras}
\label{subsec:cq}

In \cite{Nakanishi10c},
 based on  cluster algebras by \cite{Fomin02,Fomin07},
an identity of the Rogers dilogarithm was associated
 with any period of seeds of a cluster algebra.
 It looks as follows
 \begin{gather}
\label{eq:DI3int}
\sum_{t=1}^L
\varepsilon_t
L\left(
\frac{y_{k_t}(t)^{\varepsilon_t}}
{1+y_{k_t}(t)^{\varepsilon_t}}
\right)=0.
\end{gather}
A precise account will be given
in Section \ref{subsec:CDI}.
Here we only mention that
$\varepsilon_1, \dots, \varepsilon_L$
is a  certain sequence of signs
called the {\em tropical sign-sequence}.
The simplest case of the cluster algebra of type~$A_2$
yields   the pentagon relation \eqref{eq:Lpent}.
Thus, it provides a vast generalization of \eqref{eq:Lpent}.
Here we  call this family the {\em classical dilogarithm identities}.

Cluster algebras have the quantum counterparts, called
{\em quantum cluster algebras\/} \cite{Berenstein05b,Fock03}.
Here we use the formulation by \cite{Fock03}.
Any period of seeds of a classical
(nonquantum) cluster algebra
is also a period of seeds of the corresponding quantum cluster algebra
and {\em vice versa}.
Recently, in parallel with the classical case,
an identity  of the quantum dilogarithm $\mathbf{\Psi}_q(x)$
was associated with
any period of seeds of a quantum cluster algebra
by \cite{Keller11} (see also \cite{Reineke08,Nagao11,Nagao11b}).
Moreover, as a pleasant surprise,
we simultaneously obtain at least
{\em four variations\/} of quantum
dilogarithm identities as follows.

1) {\em Identities in  tropical form for $\mathbf{\Psi}_q(x)$}.
This is the form presented by \cite{Keller11},
and it looks as follows
\begin{gather}
\label{eq:QDI2int}
\mathbf{\Psi}_q(\mathsf{Y}^{\varepsilon_1 \alpha_1}
 )^{\varepsilon_1}
\cdots
\mathbf{\Psi}_q(\mathsf{Y}^{\varepsilon_L \alpha_L}
 )^{\varepsilon_L}
=1.
\end{gather}
A precise account will be given in Section \ref{subsec:QDIT}.
The simplest case of the quantum cluster algebra of type $A_2$
yields the pentagon relation \eqref{eq:pent2}.

2) {\em Identities in  universal form for $\mathbf{\Psi}_q(x)$}.
This is the form presented by \cite{Volkov11,Volkov11b},
and it looks as follows
\begin{gather}
\label{eq:QDI3int}
\mathbf{\Psi}_q(Y_{k_L}(L)^{\varepsilon_L}
 )^{\varepsilon_L}
\cdots
\mathbf{\Psi}_q(Y_{k_1}(1)^{\varepsilon_1}
 )^{\varepsilon_1}
=1.
\end{gather}
A precise account will be given in Section \ref{subsec:QDIUF}.
The simplest case of type $A_2$
yields  the new variation of the
pentagon relation for $\mathbf{\Psi}_q(x)$
recently found by  \cite{Volkov11} with a suitable identif\/ication
of variables.
In general, they are obtained from the identities in tropical form
\eqref{eq:QDI2int} by
the `shuf\/f\/le method' due to A.Yu.~Volkov \cite{Volkov11b}.

3) {\em Identities in  tropical form for $\mathbf{\Phi}_b(z)$}.
This is the counterpart of the form \eqref{eq:QDI2int},
and it looks as follows
\begin{gather}
\label{eq:QDI6int}
\mathbf{\Phi}_b(\varepsilon_1 \alpha_1 \hat {\mathsf{D}}
 )^{\varepsilon_1}
\cdots
\mathbf{\Phi}_b(\varepsilon_L \alpha_L \hat {\mathsf{D}}
 )^{\varepsilon_L}
=1.
\end{gather}
A precise account will be given in Section \ref{subsec:QDITF2}.
The simplest case of type $A_2$
yields the pentagon relation \eqref{eq:pent4}.

4) {\em Identities in  local form for $\mathbf{\Phi}_b(z)$}.
This is the form presented by \cite{Fock07,Goncharov07}.
In general they are specif\/ied not only by a period of seeds but also by
 any choice of sign-sequence.
The case of tropical sign-sequence is important for our purpose,
and in that case it looks as follows
\begin{gather}
\label{eq:QDI7int}
\mathbf{\Phi}_b(\varepsilon_{1}{\hat {\mathsf{D}}}_{k_1}(1))^{\varepsilon_1}
 \rho^*_{k_1,\varepsilon_1}
  \cdots
\mathbf{\Phi}_b(\varepsilon_{L}{\hat {\mathsf{D}}}_{k_L}(L))^{\varepsilon_{L}}
 \rho^*_{k_L,\varepsilon_L}
\nu^*
=1.
\end{gather}
A precise account will be given in Section \ref{subsec:QDIiLF}.

We  call these identities \eqref{eq:QDI2int}--\eqref{eq:QDI7int}
 together the {\em quantum dilogarithm identities}.

With these classical
 and the corresponding quantum dilogarithm
identities, it is natural to ask {\em
how the latter reduce to the former
in the limit $q \rightarrow 1$ or $b\rightarrow 0$}. In this paper we  address this question.
More precisely, we demonstrate how in the limit $b \rightarrow 0$
the classical dilogarithm identities \eqref{eq:DI3int}
emerge
as the leading term in the asymptotic expansion
 from the quantum dilogarithm identities
in the form \eqref{eq:QDI7int},
that is,
{\em the local form with tropical sign-sequence.}
To do it, we  apply the saddle point
method (see, e.g., \cite[p.~95]{Takhtajan08}), also known as the stationary phase method, \`a la \cite{Faddeev94}.
In particular, we show transparently
 how the aforementioned
logarithmic gap between the Euler and Rogers
dilogarithms is f\/illed.
See Section \ref{subsec:res} for the bottom line.

Three remarks follow.
First, the variables of quantum cluster algebras
admit a natural quantum-mechanical formulation,
where the limit $b\rightarrow 0$
corresponds to the limit $\hbar \rightarrow 0$
of the Planck constant $\hbar$.
See \eqref{eq:hbar} and \eqref{eq:wD}.
Furthermore, the classical dilogarithm identities
appear as the leading terms of the
quantum dilogarithm identities
for the asymptotic expansion in $\hbar$.
Therefore, following the standard terminology of
quantum mechanics, we call the limiting procedure
the {\em semiclassical limit}.

Second,
 even though our treatment of the saddle point method
here is standard in quantum mechanics,
we admit and stress that
we did not pursue the complete, functional-analytic rigorousness.
Namely, the validity of the method in total and specif\/ic details,
for example,
 the uniqueness of the solution of the saddle point equations,
 the specif\/ication of the
integration contour through the saddle point, etc., are not argued.
Our objective here is not to prove
the classical dilogarithm identities by this method,
but to make a direct bridge between
 the classical and the quantum dilogarithm identities.

Third, there is actually the {\em fifth} form
 of quantum dilogarithm identities,
namely,
the {\em identities in local form for $\mathbf{\Psi}_q(x)$
with tropical sign-sequence}.
This is the counterpart of the form
\eqref{eq:QDI7int},
and it looks as follows
\begin{gather}
\label{eq:QDI7aint}
\mathbf{\Psi}_q(\mathsf{\hat{Y}}_{k_1}(1)^{\varepsilon_{1}})^{\varepsilon_1}
 \rho^*_{k_1,\varepsilon_1}
  \cdots
\mathbf{\Psi}_q(\mathsf{\hat{Y}}_{k_L}(L)^{\varepsilon_{L}})^{\varepsilon_{L}}
 \rho^*_{k_L,\varepsilon_L}
\nu^*
=1.
\end{gather}
One can also obtain the classical dilogarithm identities
\eqref{eq:DI3int} from
\eqref{eq:QDI7aint}
in the semiclassical limit.
However, the relevant dif\/ferential operators
are not self-adjoint. Therefore, the semiclassical limit
is more natural for $\mathbf{\Phi}_b(z)$ {\em from the quantum-mechanical
point of view}.
For the completeness, we also present it in Appendix~\ref{appendixA}.

In summary, our result
establishes the following
scheme
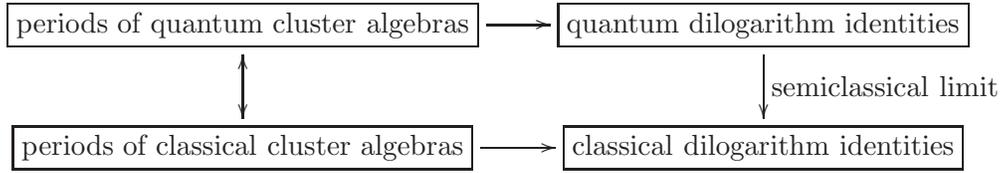
\begin{gather*}
\begin{xy}
\xymatrix{
\framebox{periods of quantum cluster algebras}
\ar@{->}[r]&
\framebox{quantum dilogarithm identities}
\ar@{->}[d]^{\mbox{semiclassical limit}}
\\
\framebox{periods of classical cluster algebras}
\ar@{->}[r] \ar@{<->}[u] &
\framebox{classical dilogarithm identities}
\\
}
\end{xy}
\end{gather*}

The organization of the paper is the following.
In Section~\ref{section2} we present the classical di\-logarithm identities
obtained from periods of cluster algebras.
In Section~\ref{section3} we present the quantum dilogarithm identities for
the quantum dilogarithm $\mathbf{\Psi}_q(x)$.
In Section~\ref{sec:phi} we present the quantum dilogarithm identities for
the Faddeev's quantum dilogarithm $\mathbf{\Phi}_b(z)$.
In Section~\ref{sec:from} we demonstrate how Rogers dilogarithm identities naturally
emerge from the quantum dilogarithm identities in local form
 in the semiclassical limit by applying the saddle point method.
This is the main part of the paper.
In Appendix~\ref{appendixA},  we present the quantum dilogarithm identities
in local form for $\mathbf{\Psi}_q(x)$.
Then, we derive the classical dilogarithm identities from them
in the semiclassical
limit.

\section{Classical dilogarithm identities}\label{section2}

In this section  we present the classical dilogarithm identities
obtained from periods of cluster algebras following \cite{Nakanishi10c}.

\subsection{Euler and Rogers dilogarithms}
\label{subsec:dilog}

Let $\mathrm{Li}_2(x)$ and $L(x)$ be the Euler and Rogers dilogarithm functions, respectively  \cite{Lewin81},
\begin{gather}
\label{eq:L00}
\mathrm{Li}_2(x)
 =-\int_{0}^x
\left\{ \frac{\log(1-y)}{y}
\right\} dy,
\qquad x\leq 1,
\\
\label{eq:L0}
L(x) =-\frac{1}{2}\int_{0}^x
\left\{ \frac{\log(1-y)}{y}+
\frac{\log y}{1-y}
\right\} dy,
\qquad  0 \leq x\leq 1.
\end{gather}
Two functions are related as follows
\begin{gather}
\label{eq:LL1}
L(x)
 =\mathrm{Li}_2(x)
+ \frac{1}{2}
\log x \log (1-x),
\qquad  0 \leq x \leq 1,
\\
\label{eq:LL2}
-L\left(\frac{x}{1+x}\right)
 =\mathrm{Li}_2(-x)
+ \frac{1}{2}
 \log x \log (1+x),
\qquad  0 \leq x.
\end{gather}
The function $L(x)$ satisf\/ies the
property \eqref{eq:Lpent} and also
the following ones
\begin{gather}
L(0)=0,
\qquad L(1)=\frac{\pi^2}{6},\nonumber\\
\label{eq:euler}
L(x)+L(1-x)=\frac{\pi^2}{6},
\qquad 0\leq x\leq 1.
\end{gather}

\subsection[$y$-variables in cluster algebras]{$\boldsymbol{y}$-variables in cluster algebras}

In this subsection
we recall some def\/initions
and properties
of the cluster algebras with
coef\/f\/i\-cients~\cite{Fomin02,Fomin03a},
following the convention of \cite{Fomin07}
with slight change of notation and terminology.
Here, we concentrate on the `coef\/f\/i\-cients' or
`$y$-variables', since we do not explicitly use the
`cluster variables' or `$x$-variables'.

Let $I$ be a f\/inite set,
and f\/ix the {\em initial $y$-seed}
$(B,y)$,
which is a pair of
a skew-symmetric (integer) matrix
$B=(b_{ij})_{i,j\in I}$
and an $I$-tuple of  commutative variables
$y=(y_i)_{i\in I}$.
Let $\mathbb{P}_{\mathrm{univ}}(y)$
be the {\em universal semifield of~$y$},
which
consists of all nonzero rational functions of~$y$
having subtraction-free expressions.
It is a semif\/ield,
i.e., the Abelian multiplicative group with addition (but not with subtraction), by the ordinary multiplication
and addition of rational functions.

Let $(B',y')$ be any pair
of
a skew-symmetric  matrix
$B'=(b'_{ij})_{i,j\in I}$
and an $I$-tuple
$y'=(y'_i)_{i\in I}$ with $y'_i\in
\mathbb{P}_{\mathrm{univ}}(y)$.
For each $k\in I$,
we def\/ine another pair $(B'',y'')=\mu_k(B',y')$
of
a skew-symmetric matrix
$B''=(b''_{ij})_{i,j\in I}$
and an $I$-tuple
$y''=(y''_i)_{i\in I}$ with $y''_i\in
\mathbb{P}_{\mathrm{univ}}(y)$,
called the {\em mutation
of $(B',y')$ at $k$},
by the following rule:

{\it  $(i)$ Mutation of matrix:}
\begin{gather}
\label{eq:Bmut}
b''_{ij}=
\begin{cases}
-b'_{ij},& \mbox{$i=k$ or $j=k$,}\\
\!\!
\begin{array}{l}
b'_{ij}+[-b'_{ik}]_+ b'_{kj}
+ b'_{ik}[b'_{kj}]_+
\\
= b'_{ij}+[b'_{ik}]_+ b'_{kj}
+ b'_{ik}[-b'_{kj}]_+,
\end{array}
&\mbox{otherwise}.
\end{cases}
\end{gather}

{\it $(ii)$  Exchange relation of $y$-variables:}
\begin{gather}
\label{eq:coef}
y''_i =
\begin{cases}
\displaystyle
{y'_k}{}^{-1},&i=k,\\
\!\!
\begin{array}{l}
y'_i y'_k{}^{[b'_{ki}]_+}(1+ y'_k)^{-b'_{ki}}
\\
=
y'_i y'_k{}^{[-b'_{ki}]_+}(1+ y'_k{}^{-1})^{-b'_{ki}},
\end{array}
&
i\neq k.\\
\end{cases}
\end{gather}
Here,
 $[a]_+=a$ for $a\geq 0$ and $0$ for $a<0$.
Starting from the initial $y$-seed
$(B,y)$, repeat the mutations.
Each  resulting pair $(B',y')$ is
 called a {\em $y$-seed of $(B,y)$}.

\begin{Remark}
The convention of \cite{Fomin07} adopted here
 is related with the  convention of
\cite{Fock03,Fock07,Keller11}
 by exchanging the matrix
$B'$ with its transposition.
\end{Remark}

\subsection[Tropical $y$-variables]{Tropical $\boldsymbol{y}$-variables}

Let $\mathbb{P}_{\mathrm{trop}}(y)$
 be the {\em tropical semifield\/}
of $y=(y_i)_{i\in I}$, which
is the Abelian multiplicative group freely generated by
$y$ endowed with the addition $\oplus$
\begin{gather*}
\prod_{i\in I} y_i^{a_i}\oplus
\prod_{i\in I} y_i^{b_i}
=
\prod_{i\in I} y_i^{\min(a_i,b_i)}.
\end{gather*}
There is a canonical surjective semif\/ield homomorphism
$\pi_{\mathbf{T}}$ (the {\em tropical evaluation})
from $\mathbb{P}_{\mathrm{univ}}(y)$
to $\mathbb{P}_{\mathrm{trop}}(y)$ def\/ined by $\pi_{\mathbf{T}}(y_i)= y_i$
and $\pi_{\mathbf{T}}(\alpha)=1$ ($\alpha \in \mathbb{Q}_+$).
For any $y$-variable~$y'_i$ of a $y$-seed
$(B',y')$ of $(B,y)$,
let us write $[y'_i]:= \pi_{\mathbf{T}}(y'_i)$ for simplicity.
We call $[y'_i]$'s the {\em tropical $y$-variables\/}
(the {\em principal coefficients\/} in \cite{Fomin07}).
They satisfy the exchange relation \eqref{eq:coef}
by replacing $y'_i$ and $+$ with $[y'_i]$
and $\oplus$.

We say that a Laurent monomial $[y'_i]$ is {\em positive\/} (resp.
 {\em negative}) if it is not 1 and all the exponents are nonnegative
(resp. nonpositive).

\begin{Proposition}[Sign-coherence \cite{Derksen10, Plamondon10b, Nagao10}]
\label{prop:sign}
For any $y$-seed $(B',y')$ of $(B,y)$,
the Laurent monomial $[y'_i]$ in $y$ is either
positive or negative.
\end{Proposition}

Based on Proposition \ref{prop:sign},
for any $y$-seed $(B',y')$ of $(B,y)$,
let $\varepsilon(y'_i)$ be 1 (resp.
$-1$) if
$[y'_i]$ is positive (resp. negative).
We call it the {\em tropical sign of $y'_i$}
by identifying $\pm 1$ with the signs $\pm$.

Using the tropical sign $\varepsilon(y'_i)$,
 the tropical exchange relation is written as follows:
\begin{gather}
\label{eq:tropcoef}
[y''_i] =
\begin{cases}
\displaystyle
[{y'_k}]^{-1},&i=k,\\
[y'_i]
 [y'_k]^{[\varepsilon(y'_k) b'_{ki}]_+},
&
i\neq k.\\
\end{cases}
\end{gather}

\subsection[Periodicity of $y$-seeds]{Periodicity of $\boldsymbol{y}$-seeds}

For any $I$-sequence $(k_1,k_2,\dots,k_L)$,
set $(B(1),y(1)):=(B,y)$,
and
consider the sequence of mutations of $y$-seeds
of $(B,y)$,
\begin{gather}
\label{eq:seedmutseq}
(B(1),y(1))
\
\mathop{\longleftrightarrow}^{\mu_{k_1}}
\
(B(2),y(2))
\
\mathop{\longleftrightarrow}^{\mu_{k_2}}
\
\cdots
\
\mathop{\longleftrightarrow}^{\mu_{k_L}}
\
(B(L+1),y(L+1)).
\end{gather}

\begin{Definition}
Let $\nu:I\rightarrow I$ be any bijection.
We say that an $I$-sequence $(k_1,k_2,\dots,k_L)$
is a~{\em $\nu$-period of $(B,y)$}
if the following holds
\begin{gather}
\label{eq:By-period}
b_{\nu(i)\nu(j)}(L+1)=b_{ij}(1),
\qquad
y_{\nu(i)}(L+1)=y_i(1),
\qquad i,j\in I.
\end{gather}
\end{Definition}

See \cite{Fomin03b,Keller10,Inoue10a,Inoue10b,Nakanishi10b}
for various examples of periodicity.

Remarkably, the periodicity of  $y$-seeds reduces to
the periodicity of tropical $y$-variables,
which is much simpler.

\begin{Proposition}[\cite{Inoue10a,Plamondon10b}]
The condition \eqref{eq:By-period} holds
if and only if
\begin{gather*}
[y_{\nu(i)}(L+1)]=[y_i(1)],
\qquad i\in I.
\end{gather*}
\end{Proposition}

For $\tilde{I}\supset I $ and a
skew-symmetric matrix $\tilde{B}=
(\tilde{b}_{ij})_{i,j\in \tilde{I}}$,
we say that
{\em $\tilde{B}$ is an $\tilde{I}$-extension of $B$}
if $\tilde{b}_{ij}=b_{ij}$ for any $i,j\in I$.

\begin{Example}
\label{ex:ext}
For any skew-symmetric matrix $B$ with index set $I$,
which may be degenerate,
let ${I}'=\{i'\mid i\in I\}$ be
a copy of $I$ and let
$\tilde{I}=I\sqcup I'$.
Def\/ine the skew-symmetric matrix $\tilde{B}=
(\tilde{b}_{ij})_{i,j\in \tilde{I}}$ by
\begin{gather*}
\tilde{b}_{ij}=
\begin{cases}
b_{ij},&i,j\in I,\\
1, & j\in I, i=j',\\
-1, & i\in I, j=i',\\
0, & \mbox{otherwise}.
\end{cases}
\end{gather*}
Then, $\tilde{B}$ is an $\tilde{I}$-extension of $B$;
furthermore, $\tilde{B}$ is {\em nondegenerate}.
The matrix $\tilde{B}$ is called the {\em principal
extension of $B$}.
\end{Example}

\begin{Proposition}[{Extension Theorem
\cite{Nakanishi10c}}]
\label{prop:ext}
Suppose that an $I$-sequence $(k_1,\dots,k_L)$
is a $\nu$-period of $(B,y)$.
Then, for any   $\tilde{I}$-extension $\tilde{B}$ of
$B$,
$(k_1,\dots,k_L)$ is also a $\nu$-period of $(\tilde{B},
\tilde{y})$.
\end{Proposition}

In Proposition \ref{prop:ext}
the periodicity
of the `external' variables
$\tilde{y}_i$ ($i\in \tilde{I}\setminus I$)
is nontrivial.

\subsection{Classical dilogarithm identities}
\label{subsec:CDI}

Let $(k_1,\dots,k_L)$ be a $\nu$-period of
$(B,y)$.
For the mutation sequence
\eqref{eq:seedmutseq},
let $N_+$ and $N_-$ be the numbers of the
positive and negative monomials
among $[y_{k_1}(1)]$, \dots, $[y_{k_L}(L)]$, respectively,
so that $N_+ + N_- = L$.

The following is a generalization of the identities
\cite{
Gliozzi95,Gliozzi96,Frenkel95,Chapoton05,
Nakanishi09,Inoue10a,Inoue10b,Nakanishi10b}
originated from the central charge identities in conformal f\/ield theory
\cite{Kirillov89,Kirillov90,Bazhanov90,Kuniba93a}.

\begin{Theorem}[Classical dilogarithm identities
\cite{Nakanishi10c}]
\label{thm:DI}
The following identities hold
\begin{gather}
\label{eq:DI}
\frac{6}{\pi^2}
\sum_{
t=1}^L
L\left(
\frac{y_{k_t}(t)}{1+y_{k_t}(t)}
\right)
 =N_-,\\
\label{eq:DI'}
\frac{6}{\pi^2}
\sum_{
t=1
}^L
L\left(
\frac{1}{1+y_{k_t}(t)}
\right)
 =N_+,
\end{gather}
where the initial variables $y_i$ $(i\in I)$
arbitrarily take values in positive real numbers.
\end{Theorem}

Two identities
\eqref{eq:DI} and \eqref{eq:DI'}
are equivalent due to \eqref{eq:euler}.

\begin{Remark}
In \cite[Theorems 4.3 \& 6.4]{Nakanishi10c},
Proposition \ref{prop:ext} and Theorem \ref{thm:DI}
are stated only for $\nu=\mathrm{id}$.
However, the proofs therein
 are  also applicable  to a general $\nu$.
\end{Remark}

We introduce the sign-sequence
 $(\varepsilon_1,
\dots, \varepsilon_L)$ so that
$\varepsilon_t$ is the tropical sign of $y_{k_t}(t)$.
We call it the {\em tropical sign-sequence} of \eqref{eq:seedmutseq}.
Using \eqref{eq:euler}, one can  also rewrite
\eqref{eq:DI} and \eqref{eq:DI'}
 in the following way.

\begin{Theorem}
\label{thm:CDI2}
 For the tropical sign-sequence
 $(\varepsilon_1,\dots, \varepsilon_L)$,
\begin{gather}
\label{eq:DI3}
\sum_{t=1}^L
\varepsilon_t
L\left(
\frac{y_{k_t}(t)^{\varepsilon_t}}
{1+y_{k_t}(t)^{\varepsilon_t}}
\right)=0.
\end{gather}
\end{Theorem}

\subsection[Example of type $A_1$]{Example of type $\boldsymbol{A_1}$}
\label{subsec:a11}
Consider the simplest case, $I=\{1\}$ and
\begin{gather*}
B=(0).
\end{gather*}
Let $(k_1, k_2) = (1,1)$,
and consider the sequence of mutations of $y$-seeds
of $(B,y)$,
\begin{gather*}
 (B(1),y(1))
\
\mathop{\longleftrightarrow}^{\mu_{1}}
\
(B(2),y(2))
\
\mathop{\longleftrightarrow}^{\mu_{1}}
\
(B(3),y(3)).
\end{gather*}
Then,
\begin{gather*}
y_1(1)=y_1,
\qquad
y_1(2)=y_1^{-1},
\qquad
y_1(3)=y_1.
\end{gather*}
Thus, $(k_1, k_2)$ is a $\nu$-period with $\nu=\mathrm{id}$,
which is nothing but the involution property of the mutation.
Also
\begin{gather*}
[y_1(1)]=y_1,
\qquad
[y_1(2)]=y_1^{-1},
\qquad
[y_1(3)]=y_1
\end{gather*}
and
\begin{gather*}
\varepsilon_1 = 1,
\qquad
\varepsilon_2 = -1.
\end{gather*}
The classical dilogarithm identity \eqref{eq:DI3} is
\begin{gather*}
L \left(
\frac{y_{1}}
{1+y_{1}}
\right)
-
L\left(
\frac{y_{1}}
{1+y_1}
\right)=0,
\end{gather*}
which is trivial.

\subsection[Example of type $A_2$]{Example of type $\boldsymbol{A_2}$}
\label{subsec:a21}
Consider the simplest nontrivial case
\begin{gather*}
B=
\begin{pmatrix}
0 & -1 \\
1 & 0\\
\end{pmatrix},
\end{gather*}
which is also represented by the quiver of type $A_2$
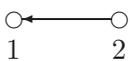
\[
\begin{picture}(40,25)(0,-15)
%
%
\put(0,0){\circle{6}}
\put(40,0){\circle{6}}
\put(37,0){\vector(-1,0){34}}
\put(-3,-15){1}
\put(37,-15){2}
\end{picture}
\]
Let $(k_1, \dots, k_5) = (1,2,1,2,1)$,
and consider the sequence of mutations of $y$-seeds
of $(B,y)$,
\begin{gather*}
 (B(1),y(1))
\
\mathop{\longleftrightarrow}^{\mu_{1}}
\
(B(2),y(2))
\
\mathop{\longleftrightarrow}^{\mu_{2}}
\
\cdots
\
\mathop{\longleftrightarrow}^{\mu_{1}}
\
(B(6),y(6)).
\end{gather*}
Then,
\begin{gather*}
\begin{cases}
y_1(1)=y_1,\\
y_2(1)=y_2,
\end{cases}
\qquad \begin{cases}
y_1(2)=y_1^{-1},\\
y_2(2)=y_2(1+y_1),
\end{cases}\qquad
\begin{cases}
y_1(3)=y_1^{-1}(1+y_2 + y_1y_2),\\
y_2(3)=y_2^{-1}(1+y_1)^{-1},
\end{cases}\nonumber\\
\begin{cases}
y_1(4)=y_1(1+y_2 + y_1y_2)^{-1},\\
y_2(4)=y_1^{-1}y_2^{-1}(1+y_2),
\end{cases}\qquad
\begin{cases}
y_1(5)=y_2^{-1},\\
y_2(5)=y_1y_2(1+y_2)^{-1},
\end{cases}\qquad
\begin{cases}
y_1(6)=y_2,\\
y_2(6)=y_1.
\end{cases}
\end{gather*}
Thus, $(k_1, \dots, k_5)$ is a $\nu$-period, where $\nu=(12)$ is
the permutation of 1 and 2.
Also
\begin{gather*}
[y_1(1)]=y_1,
\qquad
[y_2(2)]=y_2,
\qquad
[y_1(3)]=y_1^{-1},
\qquad
[y_2(4)]=y_1^{-1}y_2^{-1},
\qquad
[y_1(5)]=y_2^{-1},
\end{gather*}
and
\begin{gather*}
\varepsilon_1 = \varepsilon_2 = 1,
\qquad
\varepsilon_3 =
\varepsilon_4 =
\varepsilon_5 = -1.
\end{gather*}
The classical dilogarithm identity \eqref{eq:DI3} is
\begin{gather*}
L\left(
\frac{y_{1}}
{1+y_{1}}
\right)
+
L\left(
\frac{y_{2}(1+y_1)}
{1+y_{2}+y_1y_2}
\right)
\nonumber\\
\qquad{}-
L\left(
\frac{y_1}
{(1+y_1)(1+y_2)}
\right)
-
L\left(
\frac{y_1y_2}
{1+y_{2}+y_1y_2}
\right)
-
L\left(
\frac{y_{2}}
{1+y_{2}}
\right)
=0.
\end{gather*}
By identifying $x=y_1/(1+y_1)$, $y=y_2(1+y_1 )/(1+y_2+y_1y_2)$,
it coincides with
the pentagon relation \eqref{eq:Lpent}.

\section[Quantum dilogarithm identities for $\mathbf{\Psi}_q(x)$]{Quantum dilogarithm identities for $\boldsymbol{\mathbf{\Psi}_q(x)}$}\label{section3}

In this  section we present the quantum dilogarithm
identities for $\mathbf{\Psi}_q(x)$.
The content heavily relies on \cite{Fock03,Fock07,Keller11}.

\subsection{Quantum dilogarithm}
\label{subsec:quantum}
Following
\cite{Faddeev93,Faddeev94},
def\/ine the {\em quantum dilogarithm $\mathbf{\Psi}_q(x)$},
 for $|q|<1$ and $x\in \mathbb{C}$, by
\begin{gather}
\label{eq:psi1}
\mathbf{\Psi}_q(x) =
\sum_{n=0}^\infty
\frac{(-qx)^n}{(q^2;q^2)_n}=
\frac{1}{(-qx;q^2)_{\infty}},
\qquad
(a;q)_n = \prod_{k=0}^{n-1}\big(1-aq^{k}\big).
\end{gather}
We have the properties \eqref{eq:asym3} and \eqref{eq:pent2},
and also the following recursion relations
\begin{gather}
\label{eq:rec1}
\mathbf{\Psi}_q\big(q^{\pm2}x\big)  =
\big(1+q^{\pm 1}x\big)^{\pm1}\mathbf{\Psi}_q(x).
\end{gather}

\subsection[Quantum $y$-variables]{Quantum $\boldsymbol{y}$-variables}

So far, two kinds of
{\em quantum cluster algebras\/} are known
in the literature.
The f\/irst one
was introduced earlier by \cite{Berenstein05b},
where the $x$-variables are noncommutative
and the $y$-variables are noncommutative but restricted to
the tropical one.
The second one
was introduced by \cite{Fock03,Fock07},
where the $y$-variables are noncommutative
and the universal one
but $x$-variables are commutative.
For the relation between them,
see \cite[Section 2.7]{Fock07}
and also \cite{Tran09}.
Here  we use the second one by
\cite{Fock03,Fock07},
and concentrate on the quantum $y$-variables only.

Let $I$ be a f\/inite set, and $q$ be an
indeterminate.
We start from the {\em initial quantum $y$-seed $(B,Y)$},
which is a pair of
a skew-symmetric (integer) matrix
$B=(b_{ij})_{i,j\in I}$
and an $I$-tuple of  {\em noncommutative\/} variables
${Y}=(Y_i)_{i\in I}$ with
\begin{gather}
\label{eq:yy1}
Y_i Y_j = q^{2b_{ji}} Y_j Y_i.
\end{gather}
Accordingly, let $\mathbb{T}(B,\mathsf{Y})$ be the associated {\em quantum torus},
which is the $\mathbb{Q}(q)$-algebra
generated by the noncommutative variables
$\mathsf{Y}^{\alpha}$ ($\alpha\in \mathbb{Z}^I$) with the relations
\begin{gather*}
q^{\langle\alpha,\beta\rangle}
\mathsf{Y}^{\alpha} \mathsf{Y}^{\beta}
 = \mathsf{Y}^{\alpha+\beta},
 \qquad
\langle \alpha,\beta \rangle=
-\langle \beta,\alpha \rangle= {}^t\!\alpha B \beta.
\end{gather*}
Thus, we have $\mathsf{Y}^\alpha \mathsf{Y}^\beta
=q^{2\langle\beta,\alpha\rangle}\mathsf{Y}^\beta
  \mathsf{Y}^\alpha$.
Set $\mathsf{Y}_i
:=\mathsf{Y}^{e_i}$ for the standard
unit vector $e_i$ ($i\in I$).
Then,
by identifying $\mathsf{Y}_i$ with $Y_i$,
 we recover \eqref{eq:yy1}.

Following \cite{Keller11},
let $\mathbb{A}(B,\mathsf{Y})$ be the associated {\em quantum affine space},
which is the $\mathbb{Q}(q)$-subal\-gebra
of $\mathbb{T}(B,\mathsf{Y})$
generated by
$\mathsf{Y}^{\alpha}$'s with
$\alpha\in (\mathbb{Z}_{\geq 0})^I$.
Let $\hat{\mathbb{A}}(B,\mathsf{Y})$ be the
completion of ${\mathbb{A}}(B,\mathsf{Y})$,
which consists of the noncommutative formal power series
of $\mathsf{Y}_i$'s.
The {\em complete quantum torus
$\hat{\mathbb{T}}(B,\mathsf{Y})$}
is the localization of $\hat{\mathbb{A}}(B,\mathsf{Y})$ at
 $\mathsf{Y}^{\alpha}$'s
with $\alpha\in (\mathbb{Z}_{\geq 0})^I$.
Let $\mathsf{Frac}(\mathbb{A}(B,\mathsf{Y}))$ be the noncommutative fraction f\/ield of the algebra
$\mathbb{A}(B,\mathsf{Y})$, which is viewed
as a subskewf\/ield of $\hat{\mathbb{T}}(B,\mathsf{Y})$ \cite{Berenstein05b}.

Let $(B',Y')$ be any pair
of
a skew-symmetric  matrix
$B'=(b'_{ij})_{i,j\in I}$
and an $I$-tuple
$Y'=(Y'_i)_{i\in I}$ with $Y'_i\in
\mathsf{Frac}(\mathbb{A}(B,\mathsf{Y}))$ satisfying the relations~\eqref{eq:yy1} where everything is primed.
For each $k\in I$,
we def\/ine another same kind of
 pair $(B'',Y'')=\mu_k(B',Y')$,
 called the {\em mutation
of $(B',Y')$ at $k$},
where
$B''=(b''_{ij})_{i,j\in I}$,
which is the same as \eqref{eq:Bmut},
and $Y''=(Y''_i)_{i\in I}$, $Y''_i\in
\mathsf{Frac}(\mathbb{A}(B,\mathsf{Y}))$
is given by the following rule \cite{Fock03,Fock07}:

{\it Exchange relation of quantum $y$-variables}
\begin{gather}
\label{eq:coefq}
Y''_i =
\begin{cases}
\displaystyle
{Y'_k}{}^{-1},&i=k,\\
\begin{array}{ll}
\displaystyle
\! \! q^{
b'_{ik}[b'_{ki}]_+
}
Y'_i Y'_k{}^{[b'_{ki}]_+}
\prod_{m=1}^{|b'_{ki}|}
\big(1+ q^{-\mathrm{sgn}(b'_{ki})(2m-1)}
Y'_k\big)^{-\mathrm{sgn}({b'_{ki}})}
\\
\! \!
=
\displaystyle
q^{
b'_{ik}[-b'_{ki}]_+
}
Y'_i Y'_k{}^{[-b'_{ki}]_+}
\prod_{m=1}^{|b'_{ki}|}
\big(1+ q^{\mathrm{sgn}(b'_{ki})(2m-1)}
Y'_k{}^{-1}\big)^{-\mathrm{sgn}({b'_{ki}})},
\end{array}
&
i\neq k.
\end{cases}
\end{gather}
Formally setting $q=1$, it reduces to
\eqref{eq:coef}.

Now, starting from the quantum initial $y$-seed
$(B,Y)$, repeat the mutations.
Each  resulting pair $(B',Y')$ is
 called a {\em quantum $y$-seed of $(B,Y)$}.

\subsection{Decomposition of mutations}
\label{subsec:mut}

Let $(B',Y')$ and  $(B'',Y'')$  be  a pair of  quantum $y$-seeds of
$(B,Y)$ such that $(B'',Y'')=\mu_k(B',Y')$.
Following \cite{Fock03}, we decompose the mutation
\eqref{eq:coefq} into two parts, namely,
the monomial part and the automorphism part.

{\bf (a) Monomial part.} Def\/ine the isomorphisms $\tau_{k,\varepsilon}$
for each $\varepsilon = \pm 1$ by
\begin{gather}
\tau_{k,\varepsilon}: \
\mathsf{Frac}(\mathbb{A}(B'',\mathsf{Y}''))
 \rightarrow \mathsf{Frac}(\mathbb{A}(B',\mathsf{Y}')),
\nonumber \\
\phantom{\tau_{k,\varepsilon}: \ }{} \ \mathsf{Y}''^{}_i
 \mapsto
\begin{cases}
\mathsf{Y}'^{-1}_k, & i=k,\\
\mathsf{Y}'^{e_i + [\varepsilon b'_{ki}]_+ e_k}, & i\neq k.
\end{cases}\label{eq:phi1}
\end{gather}
The dependence of the map $\tau_{k,\varepsilon}$ on its source $(B'',Y'')$
and target $(B',Y')$
is suppressed for the notational simplicity
and should be understood in the context.
One can check that
they are indeed homomorphisms using~\eqref{eq:Bmut};
furthermore, they are isomorphisms because the inverses are given by
$\tau_{k,-\varepsilon}$ with $b'_{ki}$ being replaced by $b''_{ki}=-b'_{ki}$
in~\eqref{eq:phi1}.
Compare with the exchange relation of tropical
$y$-variables~\eqref{eq:tropcoef}.
Also note that, in ${\mathbb{A}}(B',\mathsf{Y}')$,
\begin{gather}
\label{eq:coefq2}
\mathsf{Y}'{}^{e_i+[\varepsilon b'_{ki}]_+e_k}=
q^{
b'_{ik}[\varepsilon b'_{ki}]_+
}
\mathsf{Y}'_i \mathsf{Y}'_k{}^{[\varepsilon b'_{ki}]_+}.
\end{gather}

{\bf (b) Automorphism part.} It follows from \eqref{eq:rec1}
that, for $ \mathsf{Y}'_k\in \mathbb{A}(B',\mathsf{Y}')$,
the adjoint action $\mathsf{Ad}(\mathbf{\Psi}_q(\mathsf{Y}'_k))$
is def\/ined on $\mathsf{Frac}(\mathbb{A}(B',\mathsf{Y}'))$
by
\begin{gather}
\mathsf{Ad}(\mathbf{\Psi}_q(\mathsf{Y}'_k))
(\mathsf{Y}'_i) :=
 \mathbf{\Psi}_q(\mathsf{Y}'_k) \mathsf{Y}'_i
\mathbf{\Psi}_q(\mathsf{Y}'_k)^{-1}
= \mathsf{Y}'_i
\mathbf{\Psi}_q\big(q^{-2b'_{ki}}\mathsf{Y}'_k\big)
\mathbf{\Psi}_q(\mathsf{Y}'_k)^{-1}\nonumber\\
\phantom{\mathsf{Ad}(\mathbf{\Psi}_q(\mathsf{Y}'_k)) (\mathsf{Y}'_i)}{}
=
\mathsf{Y}'_i \prod_{m=1}^{|b'_{ki}|}
\big(1+ q^{-\mathrm{sgn}(b'_{ki})(2m-1)}
\mathsf{Y}'_k \big)^{-\mathrm{sgn}({b'_{ki}})},\label{eq:ad1}
\end{gather}
and similarly,
\begin{gather}
\mathsf{Ad}\big(\mathbf{\Psi}_q(\mathsf{Y}'^{-1}_k)^{-1}\big)
(\mathsf{Y}'_i) :=
  \mathbf{\Psi}_q\big(\mathsf{Y}'^{-1}_k\big)^{-1}
 \mathsf{Y}'_i
\mathbf{\Psi}_q\big(\mathsf{Y}'^{-1}_k\big) = \mathsf{Y}'_i
\mathbf{\Psi}_q\big(q^{2b'_{ki}}\mathsf{Y}'^{-1}_k\big)^{-1}
\mathbf{\Psi}_q\big(\mathsf{Y}'^{-1}_k\big)\nonumber\\
\phantom{\mathsf{Ad}\big(\mathbf{\Psi}_q(\mathsf{Y}'^{-1}_k)^{-1}\big) (\mathsf{Y}'_i)}{}
=
\mathsf{Y}'_i \prod_{m=1}^{|b'_{ki}|}
\big(1+ q^{\mathrm{sgn}(b'_{ki})(2m-1)}
\mathsf{Y}'^{-1}_k\big)^{-\mathrm{sgn}({b'_{ki}})}.\label{eq:ad2}
\end{gather}

By combining \eqref{eq:phi1}--\eqref{eq:ad2},
we have the following
intrinsic description of the exchange rela\-tion~\eqref{eq:coefq}.

\begin{Proposition}[{\cite{Fock03,Keller11}}]
\label{prop:mutad}
We have the equality
\begin{gather}
\label{eq:adad}
 (\mathsf{Ad}(\mathbf{\Psi}_q(\mathsf{Y}'_k))
\tau_{k,+} ) (\mathsf{Y}''_i)
=
 \big(\mathsf{Ad}\big(\mathbf{\Psi}_q\big(\mathsf{Y}'^{-1}_k\big)^{-1}\big)
\tau_{k,-} \big) (\mathsf{Y}''_i),
\end{gather}
and either side of \eqref{eq:adad} coincides with
the right hand side of
 the exchange relation \eqref{eq:coefq}
 with $Y'_i$ replaced with $\mathsf{Y}'_i$.
\end{Proposition}

\begin{Remark}
In \cite{Fock03} the case  $\varepsilon=1$ was
employed as the def\/inition of the exchange relation.
The importance of the use of {\em both} the descriptions
by $\varepsilon=\pm 1$ for quantum dilogarithm identities
was clarif\/ied by~\cite{Keller11}.
We use this ref\/inement throughout the paper.
\end{Remark}

\begin{Example}
Consider the sequence of mutations of quantum $y$-seeds
of $(B,Y)$,
\begin{gather*}
 (B(1),Y(1)):=(B,Y)
\
\mathop{\longleftrightarrow}^{\mu_{k_1}}
\
(B(2),Y(2))
\
\mathop{\longleftrightarrow}^{\mu_{k_2}}
\
(B(3),Y(3)).
\end{gather*}
Then, for any sign-sequence $(\varepsilon_1,\varepsilon_2)$, we have
\begin{gather}
Y_{i}(2) =
 \big(\mathsf{Ad}(\mathbf{\Psi}_q(\mathsf{Y}_{k_1}(1)^{\varepsilon_{1}}
 )^{\varepsilon_1})
 \tau_{k_1,\varepsilon_1}\big)
 (\mathsf{Y}_{i}(2)),\nonumber\\
\label{eq:ad4}
Y_{i}(3) =
 \big(\mathsf{Ad}(\mathbf{\Psi}_q(\mathsf{Y}_{k_1}(1)^{\varepsilon_{1}}
 )^{\varepsilon_1})
 \tau_{k_1,\varepsilon_1}
 \mathsf{Ad}(\mathbf{\Psi}_q(\mathsf{Y}_{k_2}(2)^{\varepsilon_{2}}
 )^{\varepsilon_2})
 \tau_{k_2,\varepsilon_2}
 \big)
 (\mathsf{Y}_{i}(3)).
\end{gather}
\end{Example}

\subsection{Quantum dilogarithm identities in tropical form}
\label{subsec:QDIT}

For any $I$-sequence $(k_1,k_2,\dots,k_L)$,
set $(B(1),Y(1)):=(B,Y)$,
and
consider the sequence of mutations of quantum $y$-seeds
of $(B,Y)$,
\begin{gather}
\label{eq:seedmutseqq}
 (B(1),Y(1))
\
\mathop{\longleftrightarrow}^{\mu_{k_1}}
\
(B(2),Y(2))
\
\mathop{\longleftrightarrow}^{\mu_{k_2}}
\
\cdots
\
\mathop{\longleftrightarrow}^{\mu_{k_L}}
\
(B(L+1),Y(L+1)).
\end{gather}
We say that an $I$-sequence $(k_1,k_2,\dots,k_L)$
is a {\em $\nu$-period of $(B,Y)$}
if the following condition holds
for the sequence \eqref{eq:seedmutseqq}
\begin{gather}
\label{eq:By-periodq}
b_{\nu(i)\nu(j)}(L+1)=b_{ij}(1),
\qquad
Y_{\nu(i)}(L+1)=Y_i(1),
\qquad i,j\in I.
\end{gather}

The following theorem,
essentially due to \cite{Fock07},
 tells that the periodicities
of quantum $y$-seeds and (classical) $y$-seeds coincide.
\begin{Proposition}[\cite{Fock07}]
\label{prop:BB}
The condition \eqref{eq:By-period} holds for
the sequence \eqref{eq:seedmutseq}
if and only if
the condition \eqref{eq:By-periodq} holds for
the sequence \eqref{eq:seedmutseqq}.
\end{Proposition}

\begin{proof} The `if' part  immediately follows
by formally setting $q=1$ in the exchange rela\-tion~\eqref{eq:coefq} for quantum $y$-seeds.
 The 'only if' part is
 proved by \cite[Lemma~2.22]{Fock07}
 using \cite[Theorem~6.1]{Berenstein05b},
 when the matrix $B$ is nondegenerate.
When $B$ is degenerate,
thanks to
Example \ref{ex:ext} and
Proposition \ref{prop:ext},
it is reduced to the nondegenerate case.
\end{proof}

Now suppose that $(k_1,k_2,\dots,k_L)$
is a $\nu$-period of $(B,Y)$.
Due to the periodicity of \linebreak $B_{\nu(i)\nu(j)}(L+1)=
B_{ij}(1)$,
we have the isomorphism
\begin{gather*}
\mathsf{Frac}(\mathbb{A}(B(1),\mathsf{Y}(1)))
  \rightarrow
\mathsf{Frac}(\mathbb{A}(B(L+1),\mathsf{Y}(L+1))),\\
\mathsf{Y}_i(1)  \mapsto \mathsf{Y}_{\nu(i)}(L+1).
\end{gather*}
Let $\nu$ also denote this isomorphism by abusing the
notation.
For any sign-sequence
$(\varepsilon_1,\dots,\varepsilon_k)$,
the periodicity for \eqref{eq:seedmutseqq} is expressed as follows~\cite{Keller11}.
\begin{gather}
\label{eq:ad3}
 \mathsf{Ad}(\mathbf{\Psi}_q(\mathsf{Y}_{k_1}(1)^{\varepsilon_{1}}
 )^{\varepsilon_1})
 \tau_{k_1,\varepsilon_1}
  \cdots
 \mathsf{Ad}(\mathbf{\Psi}_q(\mathsf{Y}_{k_L}(L)^{\varepsilon_{L}}
 )^{\varepsilon_L})
 \tau_{k_L,\varepsilon_L}
\nu = \mathrm{id}_{
\mathsf{Frac}(
\mathbb{A}(B(1),\mathsf{Y}(1)))}.
\end{gather}

To extract the identity involving only the quantum
dilogarithm
$\mathbf{\Psi}_q(y)$,
we have to set $(\varepsilon_1$,
\dots, $ \varepsilon_L)$ to be the tropical sign-sequence
of \eqref{eq:seedmutseq}.
(We also call it the tropical sign-sequence of~\eqref{eq:seedmutseqq}.)

The following theorem is due to \cite[Theorem~5.6]{Keller11}.
The case of simply laced f\/inite type for certain periods
 was obtained by~\cite{Reineke08} with a dif\/ferent method.
See also \cite[Comments~(a), p.~5]{Nagao11}, \cite{Nagao11b} for the connection to
the Donaldson--Thomas invariants.
We include a proof because the argument therein
will be used also later.

\begin{Theorem}[Quantum dilogarithm identities in tropical form \cite{Reineke08,Keller11}]
\label{thm:QDI1}
Suppose that  $(k_1,\dots,$ $k_L)$ is a $\nu$-period
of $(B,Y)$,
and let  $(\varepsilon_1,
\dots,  \varepsilon_L)$  be the tropical sign-sequence
of~\eqref{eq:seedmutseqq}.
Let $y_i(t)$
 be the corresponding
$($classical$)$ $y$-variables in~\eqref{eq:seedmutseq},
 and let
 $\alpha_t\in \mathbb{Z}^I$ $(t=1,\dots, L)$
 be the vectors such that
 $[y_{k_t}(t)]=y^{\alpha_t}$.
 $($The vector $\alpha_t$ is called the
 $\mathbf{c}$-vector of $y_{k_t}(t)$ in
 {\rm \cite{Fomin07}.)}
 Then, the following identity holds
\begin{gather}
\label{eq:QDI2}
\mathbf{\Psi}_q(\mathsf{Y}^{\varepsilon_1 \alpha_1}
 )^{\varepsilon_1}
\cdots
\mathbf{\Psi}_q(\mathsf{Y}^{\varepsilon_L \alpha_L}
 )^{\varepsilon_L}
=1,
\end{gather}
where $\mathsf{Y}^{\varepsilon_1 \alpha_1},
\dots,\mathsf{Y}^{\varepsilon_L \alpha_L}\in
\mathbb{A}(B,\mathsf{Y})$.
\end{Theorem}

\begin{proof}
For the choice of $\varepsilon_{t}$ above,
the periodicity of  tropical $y$-variables
implies
\begin{gather}
\label{eq:ty2}
\tau_{k_1,\varepsilon_1}\cdots
\tau_{k_L,\varepsilon_L}
\nu = \mathrm{id}.
\end{gather}
Also note that  $\mathsf{Y}_{k_1}(1)^{\varepsilon_1}
=\mathsf{Y}^{\varepsilon_1\alpha_1}$ with
  $\alpha_1=e_{k_1}$ and $\varepsilon_1=1$.
Then, push out all $\tau_{k_t,\varepsilon_t}$'s to the right
in \eqref{eq:ad3} as follows
\begin{gather*}
 \mathsf{Ad}(\mathbf{\Psi}_q(\mathsf{Y}^{\varepsilon_1\alpha_1}
 )^{\varepsilon_1})
 \tau_{k_1,\varepsilon_{1}}
  \mathsf{Ad}(\mathbf{\Psi}_q(\mathsf{Y}_{k_2}(2)^{\varepsilon_2}
 )^{\varepsilon_2})
 \tau_{k_2,\varepsilon_2}
   \mathsf{Ad}(\mathbf{\Psi}_q(\mathsf{Y}_{k_3}(3)^{\varepsilon_3}
 )^{\varepsilon_3})
\cdots
 \nu
 =\mathrm{id},\\
 \mathsf{Ad}(\mathbf{\Psi}_q(\mathsf{Y}^{\varepsilon_1\alpha_1}
 )^{\varepsilon_1})
  \mathsf{Ad}(\mathbf{\Psi}_q(\mathsf{Y}^{\varepsilon_2\alpha_2}
 )^{\varepsilon_2})
 \tau_{k_1,\varepsilon_1}
 \tau_{k_2,\varepsilon_2}
    \mathsf{Ad}(\mathbf{\Psi}_q(\mathsf{Y}_{k_3}(3)^{\varepsilon_3}
 )^{\varepsilon_3})
\cdots
 \nu
 =\mathrm{id},\\
\cdots\cdots\cdots\cdots\cdots\cdots\cdots\cdots\cdots\cdots\cdots\cdots\cdots\cdots\cdots
\cdots\cdots\cdots\cdots\cdots\cdots\cdots\cdots\cdots
 \\
 \mathsf{Ad}(\mathbf{\Psi}_q(\mathsf{Y}^{\varepsilon_1\alpha_1}
 )^{\varepsilon_1})
\cdots
   \mathsf{Ad}(
   \mathbf{\Psi}_q(\mathsf{Y}^{\varepsilon_L\alpha_L}
 )^{\varepsilon_L})
  \tau_{k_1,\varepsilon_1}
  \cdots
 \tau_{k_L,\varepsilon_L}
 \nu
 =\mathrm{id}.
\end{gather*}
Thus, thanks to \eqref{eq:ty2}, we have for any $i\in I$
\begin{gather}
\label{eq:adid}
 \mathsf{Ad}(\mathbf{\Psi}_q(\mathsf{Y}^{\varepsilon_1\alpha_1}
 )^{\varepsilon_1}
\cdots
    \mathbf{\Psi}_q(\mathsf{Y}^{\varepsilon_L\alpha_L}
 )^{\varepsilon_L})
 )
 (\mathsf{Y}_{i}(1)) =\mathsf{Y}_i(1).
 \end{gather}
If $B$ is nondegenerate,
by considering the canonical form of $B$,
one can easily show that
the only  $\mathsf{Y}^{\alpha}$ which
commutes with all $\mathsf{Y}_i$'s is~1.
Therefore, \eqref{eq:adid} implies the identity \eqref{eq:QDI2}.
If  $B$ is degenerate,
again thanks to
 Example~\ref{ex:ext} and
Proposition~\ref{prop:ext},
it is reduced to the nondegenerate case.
\end{proof}

\subsection{Quantum dilogarithm identities  in universal form}
\label{subsec:QDIUF}

Let us rewrite the identity
\eqref{eq:QDI2} with genuine `nontropical', or universal, quantum
$y$-variab\-les~$Y_{k_t}(t)$.
This  generalizes the new variation of the pentagon relation~\eqref{eq:DI8} recently found by~\cite{Volkov11}
and its generalization to any simply laced f\/inite type
\cite{Volkov11b}.
To be more precise,
the pentagon relation of \cite{Volkov11} is expressed by the
{\em initial variables of the Y-system},
while our version is expressed by the {\em initial $y$-variables},
so that they have dif\/ferent expressions.
However, they coincide under a suitable identif\/ication of variables
as shown in Section~\ref{subsec:a22}.
A.Yu.~Volkov explained us how to derive his pentagon relation
and its generalization to any simply laced f\/inite type
from the tropical one  using the `shuf\/f\/le method'
in the Y-system setting~\cite{Volkov11b}.
Below we apply his shuf\/f\/le method adapted in our cluster algebraic setting.

\begin{Lemma}
Under the same assumption of
Theorem {\rm \ref{thm:QDI1}},
the following formulas hold for $t=2,\dots,L$
$($we call \eqref{eq:tnt2} the shuffle formula$)$
\begin{gather}
\label{eq:tnt1}
\mathbf{\Psi}_q(Y_{k_t}(t)^{\varepsilon_t})^{\varepsilon_t}
 =
 \mathsf{Ad}(\mathbf{\Psi}_q(\mathsf{Y}^{\varepsilon_1\alpha_1}
 )^{\varepsilon_1}
\cdots
   \mathbf{\Psi}_q(\mathsf{Y}^{\varepsilon_{t-1}\alpha_{t-1}}
 )^{\varepsilon_{t-1}})
 (\mathbf{\Psi}_q(\mathsf{Y}^{\varepsilon_t \alpha_t})^{\varepsilon_{t}}),
\\
 \label{eq:tnt2}
 \mathbf{\Psi}_q(\mathsf{Y}^{\varepsilon_1 \alpha_1}
 )^{\varepsilon_1}
 \cdots
 \mathbf{\Psi}_q(\mathsf{Y}^{\varepsilon_t \alpha_t}
 )^{\varepsilon_t}
 =
 \mathbf{\Psi}_q(Y_{k_t}(t)^{\varepsilon_t}
 )^{\varepsilon_t}
 \cdots
 \mathbf{\Psi}_q(Y_{k_1}(1)^{\varepsilon_1}
 )^{\varepsilon_1}.
\end{gather}
\end{Lemma}
\begin{proof}
Let us  prove \eqref{eq:tnt1} for $t=3$, for example.
By setting $i=k_3$ in \eqref{eq:ad4}
and
 repeating the  argument in the proof
of Theorem \ref{thm:QDI1},
we have
\begin{gather*}
Y_{k_3}(3) =
 (\mathsf{Ad}(\mathbf{\Psi}_q(\mathsf{Y}_{k_1}(1)^{\varepsilon_{1}}
 )^{\varepsilon_1})
 \tau_{k_1,\varepsilon_1}
 \mathsf{Ad}(\mathbf{\Psi}_q(\mathsf{Y}_{k_2}(2)^{\varepsilon_{2}}
 )^{\varepsilon_2})
 \tau_{k_2,\varepsilon_2}
 )
 (\mathsf{Y}_{k_3}(3))
\\
\phantom{Y_{k_3}(3)}{} =
 \mathsf{Ad}(\mathbf{\Psi}_q(\mathsf{Y}^{\varepsilon_{1}\alpha_1}
 )^{\varepsilon_1}
\mathbf{\Psi}_q(\mathsf{Y}^{\varepsilon_{2}\alpha_2}
 )^{\varepsilon_2})
 (\mathsf{Y}^{\alpha_3}).
\end{gather*}
Then, by extending the map
$\mathsf{Ad}(\mathbf{\Psi}_q(\mathsf{Y}^{\varepsilon_{t}\alpha_t}
 )^{\varepsilon_t})$ to $\hat{\mathbb{T}}(B,\mathsf{Y})$,
 we obtain \eqref{eq:tnt1} for $t=3$. The general case is similar.
Then, \eqref{eq:tnt2} follows from
 \eqref{eq:tnt1} by induction.
\end{proof}

Applying \eqref{eq:tnt2} with $t=L$ to
the identity \eqref{eq:QDI2},
we immediately obtain the universal counterpart
of  \eqref{eq:QDI2}.

\begin{Corollary}[Quantum dilogarithm identities in universal form
 (\cite{Volkov11,Volkov11b})]
\label{cor:QDI2}
Under the same assumption of
Theorem {\rm \ref{thm:QDI1}},
the following identity holds
\begin{gather}
\label{eq:QDI3}
\mathbf{\Psi}_q(Y_{k_L}(L)^{\varepsilon_L}
 )^{\varepsilon_L}
\cdots
\mathbf{\Psi}_q(Y_{k_1}(1)^{\varepsilon_1}
 )^{\varepsilon_1}
=1.
\end{gather}
\end{Corollary}

Since \eqref{eq:tnt2} actually holds irrespective with the
periodicity of the sequence
\eqref{eq:seedmutseqq},
one can say that the two identities \eqref{eq:QDI2} and \eqref{eq:QDI3}
are equivalent.

\subsection[Example of type $A_2$]{Example of type $\boldsymbol{A_2}$}
\label{subsec:a22}
We continue to use the data in
Section~\ref{subsec:a21}.
For the initial quantum $y$-seed $(B,Y)$, we have
\begin{gather*}
Y_1 Y_2 = q^{2} Y_2 Y_1.
\end{gather*}
Consider the sequence of mutations of quantum $y$-seeds
of $(B,Y)$:
\begin{gather*}
 (B(1),Y(1))
\
\mathop{\longleftrightarrow}^{\mu_{1}}
\
(B(2),Y(2))
\
\mathop{\longleftrightarrow}^{\mu_{2}}
\
\cdots
\
\mathop{\longleftrightarrow}^{\mu_{1}}
\
(B(6),Y(6)).
\end{gather*}
Then,
\begin{gather*}
\begin{cases}
Y_1(1)=Y_1,\\
Y_2(1)=Y_2,
\end{cases}\qquad
\begin{cases}
Y_1(2)=Y_1^{-1},\\
Y_2(2)=Y_2(1+q Y_1),
\end{cases}\qquad
\begin{cases}
Y_1(3)=Y_1^{-1}(1+q Y_2 + Y_1Y_2),\\
Y_2(3)=Y_2^{-1}(1+q^{-1} Y_1)^{-1},
\end{cases}\\
\begin{cases}
Y_1(4)=Y_1(1+q Y_2 + Y_1Y_2)^{-1},\\
Y_2(4)=q ^{-1}Y_1^{-1}Y_2^{-1} (1+q Y_2), \!
\end{cases}\qquad\!\!\!\!
\begin{cases}
Y_1(5)=Y_2^{-1},\\
Y_2(5)=q^{-1}Y_1Y_2(1+q^{-1} Y_2), \!
\end{cases}\qquad\!\!\!\!
\begin{cases}
Y_1(6)=Y_2,\\
Y_2(6)=Y_1. \!
\end{cases}\!\!
\end{gather*}
The quantum dilogarithm identity in tropical form \eqref{eq:QDI2}
is
\begin{gather*}
\mathbf{\Psi}_q\left(
Y_1
\right)
\mathbf{\Psi}_q\left(
Y_2
\right)
\mathbf{\Psi}_q\left(
Y_1
\right)^{-1}
\mathbf{\Psi}_q\left(
q^{-1} Y_1 Y_2
\right)^{-1}
\mathbf{\Psi}_q\left(
Y_2
\right)^{-1}
=1,
\end{gather*}
where we used $Y^{e_1+e_2}=q^{-1} Y_1 Y_2$.
It coincides with
the pentagon relation \eqref{eq:pent2}.
The quantum dilogarithm identity in universal form \eqref{eq:QDI3}
is
\begin{gather}
 \mathbf{\Psi}_q\left(
Y_2
\right)^{-1}
\mathbf{\Psi}_q\left(
q  (1+qY_2)^{-1}Y_2 Y_1
\right)^{-1}
\mathbf{\Psi}_q\left(
(1+ qY_2 + Y_1 Y_2)^{-1}Y_1
\right)^{-1}\nonumber\\
 \qquad {}\times
\mathbf{\Psi}_q\left(
Y_2(1+qY_1)
\right)
\mathbf{\Psi}_q\left(
Y_1
\right)
=1.\label{eq:DI7}
\end{gather}
Meanwhile,
the pentagon relation in \cite{Volkov11} reads, in our convention of $\mathbf{\Psi}_q$,
\begin{gather}
\mathbf{\Psi}_q\left(
X(1+qY)^{-1}
\right)^{-1}
\mathbf{\Psi}_q\left(
q  X (1+qX + qY )^{-1}Y
\right)^{-1}
\mathbf{\Psi}_q\left(
(1+ qX)^{-1}Y
\right)^{-1}\nonumber\\
\qquad{} \times
\mathbf{\Psi}_q\left(
X
\right)
\mathbf{\Psi}_q\left(
Y
\right)
=1,\label{eq:DI8}
\end{gather}
with $YX=q^2 XY$.
Two relations \eqref{eq:DI7} and \eqref{eq:DI8} coincide by identifying
$X=Y_2(1+qY_1)$, $Y=Y_1$.

\begin{Remark}
The relation \eqref{eq:DI8}
 should be compared with the quantum pentagon relation {\em at $N$th
roots of unity} \cite{Faddeev94}, where $N$th powers of the operators
are central and they enter the relation as parameters. As was remarked
by Bazhanov and Reshetikhin in~\cite{Bazhanov95}, these
parameters are related in exactly the same way as the arguments in the
classical pentagon relation; see
\cite[equation~(3.18)]{Bazhanov95}.
The quantum pentagon relation at roots
of unity plays a central role in the construction of invariants of
links in
arbitrary 3-manifolds by using the combinatorics of triangulations~\cite{Kashaev94} and in solvable 3-dimensional lattice models of
Bazhanov and Baxter~\cite{Bazhanov93} (it is called the restricted
star-triangle identity there).
\end{Remark}

\section[Quantum dilogarithm identities for $\mathbf{\Phi}_b(z)$]{Quantum dilogarithm identities for $\boldsymbol{\mathbf{\Phi}_b(z)}$}
\label{sec:phi}

In this  section we present the quantum dilogarithm
identities for $\mathbf{\Phi}_b(x)$.
The content heavily relies on \cite{Fock07,Fock07b}.

\subsection{Faddeev's quantum dilogarithm}
\label{subsec:faddeev}
Let $b$ be a complex number with nonzero real part.
Set
\begin{gather}
\label{eq:cq}
c_b = (b+b^{-1})\sqrt{-1}/2,\ \quad
q = e^{\pi b^2\sqrt{-1} }, \ \quad
q^{\vee} = e^{\pi b^{-2}\sqrt{-1} },\ \quad
\overline{q}=(q^{\vee})^{-1}=
 e^{-\pi b^{-2}\sqrt{-1} }.\!\!
\end{gather}
Following
\cite{Faddeev95,Faddeev94b},
def\/ine the {\em Faddeev's quantum dilogarithm $\mathbf{\Phi}_b(z)$}
for $z\in \mathbb{C}$ in the strip
 $|\mathrm{Im}\, z | < |\mathrm{Im}\, c_b|$ by
\begin{gather}
\label{eq:phi0}
\mathbf{\Phi}_b(z) =
\exp
\left(
-\frac{1}{4}\int_{-\infty}^{\infty}
\frac{e^{-2 zx\sqrt{-1} }}
{\sinh (xb) \sinh (x/b) }
\frac{dx}{x}
\right),
\end{gather}
where the singularity at $x=0$ is circled from above.
It is analytically continued to a meromorphic function on the
entire complex plane.
We have the properties \eqref{eq:asym4} and \eqref{eq:pent4},
and also the following ones
(see, e.g., \cite{Ruijsenaars97,Woronowicz00,Faddeev01,Volkov05} for more information).

$(i)$ Symmetries:
\begin{gather}
\label{eq:duality}
\mathbf{\Phi}_b(z)
=
\mathbf{\Phi}_{b^{-1}}(z)
=
\mathbf{\Phi}_{-b}(z).
\end{gather}

$(ii)$ Recurrence relation:{\samepage
\begin{gather}
\mathbf{\Phi}_b(z \pm b \sqrt{-1}) =
\big(1+e^{2\pi b z }q^{\pm1}\big)^{\pm1}
\mathbf{\Phi}_b(z),\nonumber\\
\mathbf{\Phi}_b(z \pm b^{-1} \sqrt{-1}) =
\big(1+e^{2\pi b^{-1 } z }{(q^{\vee})}{}^{\pm1}\big)^{\pm1}
\mathbf{\Phi}_b(z).\label{eq:rec2}
\end{gather}
}

$(iii)$ Unitarity: If $b$ is real or $|b|=1$, then
\begin{gather}
\label{eq:unitary}
|\mathbf{\Phi}_b(z)| = 1,
\qquad
z\in \mathbb{R}.
\end{gather}

$(iv)$ Relation to $\mathbf{\Psi}_q(x)$: If $\mathrm{Im}\, b^2>0$, then
\begin{gather}
\label{eq:phipsi}
\mathbf{\Phi}_b(z) =
\frac{\mathbf{\Psi}_{q}(e^{2\pi b z})}
{\mathbf{\Psi}_{\overline{q}}(e^{2\pi b^{-1} z})}.
\end{gather}
Note that, if $\mathrm{Im}\, b^2>0$, then
$|q|, |\overline{q}|<1$.

\subsection[Representation of quantum $y$-variables]{Representation of quantum $\boldsymbol{y}$-variables}
\label{subsec:qy}

Let us  recall a representation of
quantum $y$-variables as dif\/ferential operators
in \cite{Fock07,Fock07b}.
We continue to use the data~\eqref{eq:cq}.
In view of~\eqref{eq:cq},
we further set
\begin{gather}
\label{eq:hbar}
\hbar = \pi b^2,
\qquad
q=e^{\hbar \sqrt{-1}}.
\end{gather}

To any quantum  $y$-seed $(B',Y')$ of  $(B,Y)$
 we associate operators $\hat u'=(\hat u'_i)_{i\in I}$ and
 $\hat p'=(\hat p'_i)_{i\in I}$ satisfying the relations
\begin{gather}
\label{eq:ud2}
[\hat u'_i,\hat u'_j]=[\hat p'_i,
\hat p'_j]=0,\qquad [\hat p'_i,\hat u'_j]=\frac{\hbar}{\sqrt{-1}}\delta_{ij}.
\end{gather}
The algebra of operators $\hat u'$ and $\hat{p}'$ has a natural representation on the Hilbert space $L^2(\mathbb{R}^I)$:
\begin{gather}
\label{eq:lc1}
(\hat u'_if)(u')=u'_if(u'),
\qquad (\hat p'_i f)(u')=\frac{\hbar}{\sqrt{-1}}
\frac{\partial f(u')}{\partial u'_i},
\qquad u'\in\mathbb{R}^I.
\end{gather}
Using Dirac's notation $f(u')=\langle u'\vert f\rangle$, we have formally
\begin{gather*}
\langle u'\vert \hat u'_i\vert f\rangle=u'_i\langle u'\vert f\rangle,\qquad
\langle u'\vert \hat p'_i\vert
 f\rangle=\frac{\hbar}{\sqrt{-1}}
\frac{\partial}{\partial u'_i}\langle u'\vert f\rangle,
\end{gather*}
or simply
\begin{gather*}
\langle u'\vert \hat u'_i=u'_i\langle u'\vert,\qquad
\langle u'\vert \hat p'_i=
\frac{\hbar}{\sqrt{-1}}\frac{\partial}{\partial u'_i}\langle u'\vert.
\end{gather*}
The set of generalized vectors $\{\langle u'\vert\}_{u'\in\mathbb{R}^I}$ will be called the \emph{local coordinates} of $(B',Y')$.

Def\/ine
\begin{gather}
\label{eq:wD}
\hat w'_i =
\sum_{j\in I} b'_{ji} \hat u'_j,
\qquad
\hat D'_i =
 \hat p'_i
 + \hat w'_i,
\qquad
{\mathsf{\hat{Y}}}'_i=\exp \hat D'_i.
\end{gather}
The following relations hold
\begin{gather}
\label{eq:dd1}
[\hat D'_i,\hat D'_j] = 2\hbar {\sqrt{-1}}b'_{ji},
\qquad
\mathsf{\hat{Y}}'_i\mathsf{\hat{Y}}'_j
=q^{2b'_{ji}}\mathsf{\hat{Y}}'_j\mathsf{\hat{Y}}'_i.
\end{gather}
Also recall the following general fact,
which is a special case of the
Baker--Campbell--Hausdorf\/f formula:
 For any noncommutative
variables $A$ and  $B$ such that $[A,B]=C$ and
$[C,A]=[C,B]=0$, we have
\begin{gather*}
e^A e^B = e^{C/2} e^{A+B}.
\end{gather*}
Thus, we have a representation of
$\mathbb{T}(B,\mathsf{Y})$ on
$L^2(\mathbb{R}^I)$
with
\begin{gather}
\label{eq:Yrep1}
\mathsf{Y}'{}^{\alpha}
\mapsto
\mathsf{\hat{Y}}'{}^{\alpha}
:=
\exp \big(
\alpha \hat D'
\big),
\qquad
\alpha \hat D' :=
\sum_{i\in I}
\alpha_i \hat D'_i.
\end{gather}

\subsection{Decomposition of mutations}
\label{subsec:mut2}
Here we present a result which is analogous to that of Section \ref{subsec:mut}.
Let $(B',Y')$ and $(B'',Y'')$ be a pair of quantum
$y$-seeds of $(B,Y)$ such that
$(B'',Y'')=\mu_k (B',Y')$.

{\bf (a) Monomial part.}
For each  $\varepsilon=\pm 1$,
consider the following map
\begin{gather}
\rho_{k,\varepsilon}: \ \mathbb{R}^{I}   \rightarrow \mathbb{R}^{I}, \nonumber\\
\phantom{\rho_{k,\varepsilon}: \ }{} \ (u')   \mapsto (u''),\nonumber\\
\phantom{\rho_{k,\varepsilon}: \ }{} \ u''_i =
\begin{cases}
-u'_k + \sum\limits_{j\in I} [-\varepsilon b'_{jk}]_+ u'_j,
& i = k,\\
u'_i, & i\neq k.
\end{cases}\label{eq:utrans}
\end{gather}
Let $\rho^*_{k,\varepsilon}$ be the induced map in the space of functions $L^2(\mathbb{R}^{I})$,
\begin{gather*}
\rho^*_{k,\varepsilon}: \
L^2(\mathbb{R}^{I})  \rightarrow
L^2(\mathbb{R}^{I}),\nonumber\\
\phantom{\rho^*_{k,\varepsilon}: \  }{} \ f  \mapsto f \circ \rho_{k,\varepsilon},
\end{gather*}
or, formally,
\begin{gather*}
\langle u'\vert \rho^*_{k,\varepsilon}=\langle \rho_{k,\varepsilon}(u')\vert=\langle u''\vert,
\end{gather*}
by which we relate the local coordinates of $(B',Y')$
and $(B'',Y'')$.

For any linear operator $\hat O$ acting on $L^2( \mathbb{R}^{I})$,
let
\begin{gather*}
\mathsf{Ad}(
\rho^*_{k,\varepsilon}
)(\hat O)
:=
\rho^*_{k,\varepsilon} \hat O
(\rho^*_{k,\varepsilon})^{-1}.
\end{gather*}
In other words, it is
def\/ined by the commutative diagram
\begin{gather*}
\begin{CD}
L^2( \mathbb{R}^{I})
@>{\rho^*_{k,\varepsilon}}>>
L^2( \mathbb{R}^{I})\\
@V{\hat O}VV @VV{\mathsf{Ad}(
\rho^*_{k,\varepsilon}
)(\hat O)}V\\
L^2( \mathbb{R}^{I})
@>{\rho^*_{k,\varepsilon}}>>
L^2( \mathbb{R}^{I}).
\end{CD}
\end{gather*}
Then, we have
\begin{gather}
\label{eq:utrans2}
\mathsf{Ad}(
\rho^*_{k,\varepsilon}
)
(\hat{u}''_i) =
\begin{cases}
-\hat{u}'_k + \sum\limits _{j\in I} [-\varepsilon b'_{jk}]_+ \hat{u}'_j,
& i = k,\\
\hat{u}'_i, & i\neq k,
\end{cases}
\\
\label{eq:wtrans2}
\mathsf{Ad}(
\rho^*_{k,\varepsilon}
)
(\hat{w}''_i)=
\begin{cases}
-\hat{w}'_k, & i= k,\\
\hat{w}'_i + [\varepsilon b'_{ki}]_+ \hat{w}'_k,
& i \neq k,
\end{cases}
\\
\label{eq:dtrans}
\mathsf{Ad}(
\rho^*_{k,\varepsilon}
)
(\hat p''_i) =
\begin{cases}
-\hat p'_k, & i= k,\\
\hat p'_i + [\varepsilon b'_{ki}]_+
\hat p'_k,
& i \neq k,
\end{cases}\\
\label{eq:Dtrans}
\mathsf{Ad}(
\rho^*_{k,\varepsilon}
)
(\hat D''_i) =
\begin{cases}
-\hat D'_k, & i= k,\\
\hat D'_i + [\varepsilon b'_{ki}]_+ \hat D'_k,
& i \neq k,
\end{cases}
\end{gather}
where \eqref{eq:Dtrans} follows from
\eqref{eq:dtrans} and
\eqref{eq:wtrans2}.
It follows from \eqref{eq:Dtrans} that
\begin{gather*}
\mathsf{Ad}(
\rho^*_{k,\varepsilon}
)
(\mathsf{\hat{Y}}''_i) =
\begin{cases}
\mathsf{\hat{Y}}'_k{}^{-1}, & i= k,\\
\mathsf{\hat{Y}}'^{e_i +  [\varepsilon b'_{ki}]_+ e_k},
& i \neq k,
\end{cases}
\end{gather*}
which coincides with~\eqref{eq:phi1}.

\begin{Remark}
\label{rem:dual}
The transformation of \eqref{eq:utrans}
is the one for the {\em $g$-vectors\/} in~\cite{Fomin07}
if $\varepsilon$ is the tropical sign of $y'_k$.
Similarly, for
$w'_i = \sum_{j\in I} b'_{ji} u'_j$,
the induced transformation
\begin{gather*}
{w}''_i =
\begin{cases}
-{w}'_k, & i= k,\\
{w}'_i + [\varepsilon b'_{ki}]_+ {w}'_k,
& i \neq k,
\end{cases}
\end{gather*}
is the one for the {\em $c$-vectors\/} in~\cite{Fomin07},
and it is the logarithmic form of
the tropical exchange relation~\eqref{eq:tropcoef}.
They
are known to be dual in the
following sense~\cite{Fock03,Nakanishi11a}
\begin{gather*}
\sum_{i\in I} u''_i w''_i = \sum_{i\in I} u'_i w'_i.
\end{gather*}
\end{Remark}

{\bf (b) Automorphism part.}
We set
\begin{gather}
\label{eq:wD3}
\hat {\mathsf{D}}'_i
= \frac{1}{2\pi b} \hat D'_i.
\end{gather}
Then, we have
\begin{gather*}
\hat{\mathsf{Y}}'_i{}^{-1}
\mathbf{\Phi}_b
(\varepsilon
\hat {\mathsf{D}}'_j
)
\hat{\mathsf{Y}}'_i
=
\mathbf{\Phi}_b
\big(
\varepsilon
\hat {\mathsf{D}}'_j
- \varepsilon \sqrt{-1} b b'_{ji}
\big).
\end{gather*}
Thus,
thanks to the recurrence relation \eqref{eq:rec2},
one obtains, for each $\varepsilon=\pm 1$,
\begin{gather*}
\mathsf{Ad}(\mathbf{\Phi}_b
({\varepsilon  \hat {\mathsf{D}}'_k}
)^{\varepsilon})
(\mathsf{\hat{Y}}'_i)
=
\mathsf{\hat{Y}}'_i \prod_{m=1}^{|b'_{ki}|}
\big(1+ q^{-\varepsilon \mathrm{sgn}(b'_{ki})(2m-1)}
\mathsf{\hat{Y}}'_k{}^{\varepsilon }\big)^{-\varepsilon \mathrm{sgn}({b'_{ki}})}
\end{gather*}
by an analogous calculation to
\eqref{eq:ad1} and \eqref{eq:ad2}.

In summary, we have a parallel statement to
Proposition \ref{prop:mutad}.

\begin{Proposition}[\cite{Fock07,Fock07b}]
\label{prop:mutad2}
We have the equality
\begin{gather}
\label{eq:adad2}
 (\mathsf{Ad}(\mathbf{\Phi}_b(\hat {\mathsf{D}}'_k))
 \mathsf{Ad}(\rho^*_{k,+} )) (\mathsf{\hat{Y}}''_i)
=
 \big(\mathsf{Ad}\big(\mathbf{\Phi}_b(-\hat {\mathsf{D}}'_k)^{-1}\big)
 \mathsf{Ad}(\rho^*_{k,-} )\big) (\mathsf{\hat{Y}}''_i),
\end{gather}
and either side of \eqref{eq:adad2} coincides with
the right hand side of
 the exchange relation \eqref{eq:coefq}
 with~$Y'_i$ replaced with~$\mathsf{\hat{Y}}'_i$.
\end{Proposition}

\subsection{Dual operators}

Following \cite{Faddeev95} and \cite{Fock07},
we def\/ine the operators
$\mathsf{\hat{Z}}'_i$
which are `dual' to $\mathsf{\hat{Y}}'_i$
in the sense of the f\/irst equality of \eqref{eq:duality}.
In the situation in \eqref{eq:wD}, we def\/ine
\begin{gather*}
{\mathsf{\hat{Z}}}{}'_i=\exp(b^{-2} \hat D'_i).
\end{gather*}
Then, the following relations hold
\begin{gather}
\mathsf{\hat{Z}}'_i\mathsf{\hat{Z}}'_j
 =(q^{\vee})^{2b'_{ji}}\mathsf{\hat{Z}}'_j\mathsf{\hat{Z}}'_i,\nonumber\\
\label{eq:YZ}
\mathsf{\hat{Y}}'_i\mathsf{\hat{Z}}'_j
 =
\mathsf{\hat{Z}}'_j\mathsf{\hat{Y}}'_i.
\end{gather}

\begin{Remark}
The duality between $\mathsf{\hat{Y}}'_i$
and $\mathsf{\hat{Z}}'_i$ is not manifest
because of our preference for $b$ over  $b^{-1}$ in
\eqref{eq:wD} through \eqref{eq:hbar}.
To see it manifestly,
we set
\begin{gather*}
\hat {\mathcal{D}}'_i =\frac{1}{\gamma}
\left(
\frac{\gamma^2}{4\pi \sqrt{-1}}
 \hat \partial'_i
 + \hat w'_i
\right),
\qquad
{\mathsf{\hat{Y}}}'_i=\exp\big(2\pi b \hat {\mathcal{D}}'_i\big),
\qquad
{\mathsf{\hat{Z}}}'_i=\exp\big(2\pi b^{-1} \hat {\mathcal{D}}'_i\big),
\end{gather*}
where $\gamma$ is an arbitrary nonzero real number
and
$ \hat \partial'_i$ satisfy
$[ \hat \partial'_i,  \hat \partial'_j]=0$
and  $[\hat \partial'_i,  \hat u'_j]=\delta_{ij}$.
The following relations hold irrespective of $\gamma$
\begin{gather*}
[\hat {\mathcal{D}}'_i,\hat {\mathcal{D}}'_j] =
 \frac{\sqrt{-1}}{2\pi }b'_{ji},
\qquad
\mathsf{\hat{Y}}'_i\mathsf{\hat{Y}}'_j
=(q)^{2b'_{ji}}\mathsf{\hat{Y}}'_j\mathsf{\hat{Y}}'_i,
\qquad
\mathsf{\hat{Z}}'_i\mathsf{\hat{Z}}'_j
=(q^{\vee})^{2b'_{ji}}\mathsf{\hat{Z}}'_j\mathsf{\hat{Z}}'_i,
\qquad
\mathsf{\hat{Y}}'_i\mathsf{\hat{Z}}'_j
=\mathsf{\hat{Z}}'_j\mathsf{\hat{Y}}'_i.
\end{gather*}
Now the duality $b \leftrightarrow b^{-1}$ is manifest.
Further setting $\gamma = 2 \pi b$,
we have $\hat {\mathcal{D}'_i} =
\hat {\mathsf{D}'_i}$ and
we recover the
operators
 $\mathsf{\hat{Y}}'_i$
and $\mathsf{\hat{Z}}'_i$
in the main text.
\end{Remark}

Due to the symmetry $b\leftrightarrow b^{-1}$ in \eqref{eq:duality}
and the above remark,
we immediately obtain the following from Proposition
\ref{prop:mutad2}.

\begin{Proposition}
\label{prop:mutad3}
We have the equality
\begin{gather}
\label{eq:adad3}
 (\mathsf{Ad}(\mathbf{\Phi}_b(\hat {\mathsf{D}}'_k))
 \mathsf{Ad}(\rho^*_{k,+} )) (\mathsf{\hat{Z}}''_i)
=
 \big(\mathsf{Ad}\big(\mathbf{\Phi}_b(-\hat {\mathsf{D}}'_k)^{-1}\big)
 \mathsf{Ad}(\rho^*_{k,-} )\big) (\mathsf{\hat{Z}}''_i),
\end{gather}
and either side of~\eqref{eq:adad3} coincides with
the right hand side of
 the exchange relation~\eqref{eq:coefq}
 with~$Y'_i$ and~$q$ replaced with $\mathsf{\hat{Z}}'_i$ and $q^{\vee}$,
respectively.
\end{Proposition}

\subsection{Quantum dilogarithm identities in tropical  form}
\label{subsec:QDITF2}
Suppose that $(k_1,k_2,\dots,k_L)$
is a $\nu$-period of $(B,Y)$ as in
Section \ref{subsec:QDIT}.

The identities parallel to
\eqref{eq:QDI2} are available for $\mathbf{\Phi}_b(z)$
 directly from \eqref{eq:QDI2}.

\begin{Theorem}[Quantum dilogarithm identities in tropical form]
\label{thm:QDI4}
Under the same assumption of Theorem {\rm \ref{thm:QDI1}}
$($in particular, $(\varepsilon_1,\dots,\varepsilon_L)$ is the tropical
sign-sequence of the mutation sequence$)$, the following identity holds.
\begin{gather}
\label{eq:QDI6}
\mathbf{\Phi}_b(\varepsilon_1 \alpha_1 \hat {\mathsf{D}}
 )^{\varepsilon_1}
\cdots
\mathbf{\Phi}_b(\varepsilon_L \alpha_L \hat {\mathsf{D}}
 )^{\varepsilon_L}
=1,
\end{gather}
where $\alpha_t \hat {\mathsf{D}}=\sum_{i\in I} \alpha_i(t) \hat
 {\mathsf{D}}_i$
and  $\hat {\mathsf{D}}_i$ is the operator in \eqref{eq:wD3} for  $(B,Y)$.
\end{Theorem}
\begin{proof}
Due to the symmetry $b\leftrightarrow b^{-1}$
in \eqref{eq:duality},
one can assume that $\mathrm{Im}\, b^2 \geq 0$
without losing generality.
By the continuity of
 $\mathbf{\Phi}_b$
with respect to $b$,
it is enough to show the claim for
 $\mathrm{Im}\, b^2 >0$.
Then, by \eqref{eq:phipsi}, we have
\begin{gather*}
 \mathbf{\Phi}_b(\varepsilon_t\alpha_t \hat {\mathsf{D}})
=
\frac{ \mathbf{\Psi}_q(\hat{\mathsf{Y}}^{\varepsilon_t\alpha_t})}
{ \mathbf{\Psi}_{\overline{q}}(\hat{\mathsf{Z}}^{\varepsilon_t\alpha_t})}.
\end{gather*}
Then, thanks to the commutativity \eqref{eq:YZ},
the relation \eqref{eq:QDI6} factorizes into
 two identities
\begin{gather}
\label{eq:QDI9}
\mathbf{\Psi}_q(\mathsf{\hat Y}^{\varepsilon_1 \alpha_1}
 )^{\varepsilon_1}
\cdots
\mathbf{\Psi}_q(\mathsf{\hat Y}^{\varepsilon_L \alpha_L}
 )^{\varepsilon_L}
=1,\\
\label{eq:QDI10}
\mathbf{\Psi}_{\overline{q}}(\mathsf{\hat Z}^{\varepsilon_L \alpha_L}
 )^{\varepsilon_L}
\cdots
\mathbf{\Psi}_{\overline{q}}(\mathsf{\hat Z}^{\varepsilon_1 \alpha_1}
 )^{\varepsilon_1}
=1,
\end{gather}
where  \eqref{eq:QDI9} is
a specialization of \eqref{eq:QDI2},
 while  \eqref{eq:QDI10} is equivalent to
\begin{gather*}
\mathbf{\Psi}_{q^{\vee}}(\mathsf{\hat Z}^{\varepsilon_1 \alpha_1}
 )^{\varepsilon_1}
\cdots
\mathbf{\Psi}_{q^{\vee}}(\mathsf{\hat Z}^{\varepsilon_L \alpha_L}
 )^{\varepsilon_L}
=1,
\end{gather*}
which is another specialization of \eqref{eq:QDI2}.
\end{proof}

\subsection{Quantum dilogarithm identities in local form}
\label{subsec:QDIiLF}

Let $\langle u(t)\vert$ and  $\hat D_i(t)$ denote
the local coordinates and the operator
in \eqref{eq:wD} for $(B(t),Y(t))$, respectively.
Let $L^2(\mathbb{R}^I)_t$ be the Hilbert space together
 with the local coordinate
$\langle u(t)\vert$,
so that $\rho^*_{k_t,\varepsilon_t}:
L^2(\mathbb{R}^I)_{t+1}\rightarrow L^2(\mathbb{R}^I)_{t}$.

For the bijection $\nu$, we apply the same formalism as
$\rho$.
Namely, let $\nu: \mathbb{R}^I \rightarrow \mathbb{R}^I$
be the coordinate transformation def\/ined by $(u(L+1))\mapsto (u(1))$ with $u_i(1)=u_{\nu(i)}(L+1)$.
Def\/ine
$\nu^*:
L^2( \mathbb{R}^{I})_1 \rightarrow
L^2( \mathbb{R}^{I})_{L+1}$,
$f \mapsto f \circ \nu$,
and $\mathsf{Ad}(\nu^*)(\hat O):= \nu^* \hat O (\nu^*)^{-1}$
 for any linear operator $\hat O$ acting on
$L^2( \mathbb{R}^{I})$.
Then, $\mathsf{Ad}(\nu^*)(\hat D_i(1))=\hat D_{\nu(i)}(L+1)$.

Let us recall the result of \cite[Theorem 5.4]{Fock07}.
{\em Suppose that $b$ is a nonzero real number.}
Note that,
$\mathbf{\Phi}_b(\varepsilon_t \hat {\mathsf{D}}'_{k_t}(t))$
is a {\em unitary\/} operator by  \eqref{eq:unitary}.
By the periodicity assumption and
Propositions~\ref{prop:mutad2} and~\ref{prop:mutad3}, we have the following equalities
for any sign-sequence $\vec{\varepsilon}=
(\varepsilon_{1}$, \dots, $\varepsilon_{L})$
\begin{gather*}
 \mathsf{Ad}(\mathbf{\Phi}_b(\varepsilon_{1}{\hat {\mathsf{D}}
}_{k_1}(1))^{\varepsilon_1}
 \rho^*_{k_1,\varepsilon_1}
  \cdots
\mathbf{\Phi}_b(\varepsilon_{L}{\hat {\mathsf{D}}
}_{k_L}(L))^{\varepsilon_{L}}
 \rho^*_{k_L,\varepsilon_L}
\nu^*) (\mathsf{\hat{Y}_i}(1))  = \mathsf{\hat{Y}_i(1)},\\
 \mathsf{Ad}(\mathbf{\Phi}_b(\varepsilon_{1}{\hat {\mathsf{D}}
}_{k_1}(1))^{\varepsilon_1}
 \rho^{*}_{k_1,\varepsilon_1}
  \cdots
\mathbf{\Phi}_b(\varepsilon_{L}{\hat {\mathsf{D}}
}_{k_L}(L))^{\varepsilon_{L}}
 \rho^*_{k_L,\varepsilon_L}
{\nu^*}) (\mathsf{\hat{Z}_i}(1))  = \mathsf{\hat{Z}_i(1)}.
\end{gather*}
This is equivalent to saying that the operator
\begin{gather*}
\hat O_{\vec{\varepsilon},b}=
\mathbf{\Phi}_b(\varepsilon_{1}{\hat {\mathsf{D}}
}_{k_1}(1))^{\varepsilon_1}
 \rho^*_{k_1,\varepsilon_1}
  \cdots
\mathbf{\Phi}_b(\varepsilon_{L}{\hat {\mathsf{D}}
}_{k_L}(L))^{\varepsilon_{L}}
 \rho^*_{k_L,\varepsilon_L}
\nu^*
\end{gather*}
commutes with $\mathsf{\hat{Y}_i(1)}$ and $\mathsf{\hat{Z}_i(1)}$
for any $i\in I$.
It was shown in \cite{Fock07} that,
when $b^2$ is irrational,
such $\hat O_{\vec{\varepsilon},b}$ is the identity operator up to
a complex scalar multiple $\lambda_{\vec{\varepsilon},b}$
by generalizing  the result of \cite{Faddeev95};
furthermore, the claim holds for rational $b^2$ as well  by continuity.
Since $\hat O_{\vec{\varepsilon},b}$ is  unitary,
we have $|\lambda_{\vec{\varepsilon},b}|=1$.
Therefore, one obtains the following {\em local form\/} of the quantum
dilogarithm identities.
We call it  so,
since it is described by the family of local coordinates
$\langle u(1)\vert, \dots, \langle u(L)\vert$ associated with the mutation sequence.

\begin{Theorem}[{Quantum dilogarithm identities in
local form \cite{Fock07}}]
\label{thm:QDI3}
 Let $b$ be  a nonzero real number.
For any
sign-sequence $
{\vec{\varepsilon}}=
(\varepsilon_{1}, \dots, \varepsilon_{L})$,
the following identity holds.
\begin{gather}
\label{eq:QDI5}
\mathbf{\Phi}_b(\varepsilon_{1}{\hat {\mathsf{D}}}_{k_1}(1))^{\varepsilon_1}
 \rho^*_{k_1,\varepsilon_1}
  \cdots
\mathbf{\Phi}_b(\varepsilon_{L}{\hat {\mathsf{D}}}_{k_L}(L))^{\varepsilon_{L}}
 \rho^*_{k_L,\varepsilon_L}
\nu^*
= \lambda_{\vec{\varepsilon},b},
\qquad |\lambda_{\vec{\varepsilon},b}|=1.
\end{gather}
\end{Theorem}

For
the tropical sign-sequence,
we have a stronger version of Theorem \ref{thm:QDI3}.
One can obtain it as a direct corollary of Theorem \ref{thm:QDI4},
and not via Theorem \ref{thm:QDI3}.
So the assumption that $b$ is real is not necessary here.
{\em This is the identity we use to  derive the corresponding
classical dilogarithm
identity.}

\begin{Theorem}
\label{thm:QDI5}
For the tropical
sign-sequence $
{\vec{\varepsilon}}=
(\varepsilon_{1}, \dots, \varepsilon_{L})$,
the following identity holds
\begin{gather}
\label{eq:QDI7}
\mathbf{\Phi}_b(\varepsilon_{1}{\hat {\mathsf{D}}}_{k_1}(1))^{\varepsilon_1}
 \rho^*_{k_1,\varepsilon_1}
  \cdots
\mathbf{\Phi}_b(\varepsilon_{L}{\hat {\mathsf{D}}}_{k_L}(L))^{\varepsilon_{L}}
 \rho^*_{k_L,\varepsilon_L}
\nu^*
=1.
\end{gather}
In particular,  $\lambda_{\vec{\varepsilon},b}=1$ for the tropical
sign sequence.
\end{Theorem}
\begin{proof}
By the duality in Remark \ref{rem:dual}, the periodicity of tropical
$y$-variables  \eqref{eq:ty2} is equivalent~to
\begin{gather*}
\rho^*_{k_1,\varepsilon_1}\cdots
\rho^*_{k_L,\varepsilon_L}
\nu^* = \mathrm{id}.
\end{gather*}
Multiply it from the right of
\eqref{eq:QDI6}.
Then, repeat the argument  in the proof of Theorem~\ref{thm:QDI1} in the inverse way.
\end{proof}

In summary, for the tropical sign-sequence
we have four forms of quantum dilogarithm identities
\eqref{eq:QDI2}, \eqref{eq:QDI3}, \eqref{eq:QDI5}, and \eqref{eq:QDI7}.
The f\/irst three identities are obtained from each other
  without referring to the seed periodicity of
 \eqref{eq:seedmutseqq}.
The last one is obtained from the rest
by assuming the tropical periodicity \eqref{eq:ty2}.

\subsection[Example of type $A_2$]{Example of type $\boldsymbol{A_2}$}
\label{subsec:a23}

We continue to use the data in
Sections \ref{subsec:a21} and \ref{subsec:a22}.

The quantum dilogarithm identity in tropical form \eqref{eq:QDI6}
is
\begin{gather*}
\mathbf{\Phi}_b(
\hat {\mathsf{D}}_1
)
\mathbf{\Phi}_b(
\hat {\mathsf{D}}_2
)
\mathbf{\Phi}_b(
\hat {\mathsf{D}}_1
)^{-1}
\mathbf{\Phi}_b(
\hat {\mathsf{D}}_1 + \hat {\mathsf{D}}_2
)^{-1}
\mathbf{\Phi}_b(
\hat {\mathsf{D}}_2
)^{-1}
=1,\\
[\hat {\mathsf{D}}_1,\hat {\mathsf{D}}_2]=\frac{\sqrt{-1}}{2\pi}.
\end{gather*}
By identifying $\hat Q=\hat {\mathsf{D}}_1$,
 $\hat P=\hat {\mathsf{D}}_2$, it coincides with
the pentagon relation \eqref{eq:pent4}.

Let us also write the relevant data for
the identity \eqref{eq:QDI7} explicitly
\begin{gather*}
\hat {\mathsf{D}}_{k_t}(t) =
\begin{cases}
\displaystyle
 \frac{1}{2\pi b}
\left(
\hat p_1(t)
 + \hat u_2(t)
\right),  & t=1,3,5,\vspace{1mm}\\
\displaystyle
 \frac{1}{2\pi b}
\left(
 \hat p_2(t)
 + \hat u_1(t)
\right),  & t=2,4.
\end{cases}
\end{gather*}
The images of $\hat u_1(t+1)$, $\hat u_2(t+1)$,
$\hat w_1(t+1)$, $\hat w_2(t+1)$ by the map
$\mathsf{Ad}(
\rho^*_{k_t,\varepsilon_t})$ are given in the order
\begin{gather*}
\begin{array}{@{}lllll}
t=1:
& -\hat u_1(1),& \hat u_2(1),& -\hat w_1(1),&  \hat w_2(1),
\\
t=2:
& \hat u_1(2),& -\hat u_2(2),& \hat w_1(2), &  -\hat w_2(2),
\\
t=3:
& -\hat u_1(3)+\hat u_2(3),& \hat u_2(3),& -\hat w_1(3),&  \hat w_2(3)+\hat w_1(3),
\\
t=4:
& \hat u_1(4),& -\hat u_2(4)+\hat u_1(4),& \hat w_1(4)+\hat w_2(4),&  -\hat w_2(4),\\
t=5:
& -\hat u_1(5)+\hat u_2(5),& \hat u_2(5),& -\hat w_1(5),&  \hat w_2(5)+\hat w_1(5).
\end{array}
\end{gather*}

\section{From quantum to classical dilogarithm identities}
\label{sec:from}

In this section we demonstrate how
the classical quantum dilogarithm identities
\eqref{eq:DI3}
emerge from the quantum dilogarithm identities
in local form  \eqref{eq:QDI7}
in the semiclassical limit.
This  is the main part of the paper.

\subsection{Position and momentum bases}

We are going to evaluate the operator
in the left hand side of \eqref{eq:QDI7},
which is actually the identity operator,
by the standard quantum physics method.

{\it Throughout Section {\rm \ref{sec:from}} we assume that
 $b$ is a nonzero real number.}

Recall that we set $\hbar=\pi b^2$ in \eqref{eq:hbar}.
The asymptotic property \eqref{eq:asym4} is written as
\begin{gather}
\label{eq:asym2}
\mathbf{\Phi}_b\left(\frac{z}{2\pi b}\right) \sim
\exp \left( \frac{\sqrt{-1}}{\hbar}\frac{1}{2}\mathrm{Li}_2(-e^{z})
\right),
\qquad \hbar \rightarrow 0,
\end{gather}
where and in the rest
 $\sim$ means the leading term for the
asymptotic expansion in $\hbar$.

Let $(B(t),Y(t))$ be the quantum $Y$-seed of $(B(1),Y(1))=(B,Y)$
in \eqref{eq:seedmutseqq}.
Let $L^2(\mathbb{R}^I)_t$ be the Hilbert space together with the local coordinate
$\langle u(t)\vert$ associated with $(B(t),Y(t))$
 in the previous section.

Let $\{ |u(t) \rangle
\mid u(t)\in \mathbb{R}^I \}$
 and $\{ |p(t) \rangle
 \mid p(t)\in \mathbb{R}^I \} $ be the standard
position and the momentum
bases of $L^2(\mathbb{R}^I)_t$, respectively.
They satisfy the following properties, where $n= |I|$,
\begin{gather}
 \hat{u}_i(t)| u(t) \rangle   =
 u_i(t) | u(t) \rangle ,
\qquad
 \hat p_i | p(t) \rangle  =
 p_i(t) | p (t) \rangle,\nonumber\\
\langle u(t)|u'(t) \rangle = \prod_{i\in I} \delta(u_i(t) - u'_i(t)),
\qquad
\langle p(t)|p'(t) \rangle = (2\pi \hbar)^n \prod_{i\in I} \delta(p_i(t) - p'_i(t)),\nonumber\\
\langle u(t)|p(t) \rangle = \exp\left( \frac{\sqrt{-1}}{\hbar} u(t) p(t)\right),
\qquad
\langle p(t)|u(t) \rangle = \exp\left( -\frac{\sqrt{-1}}{\hbar} u(t) p(t)\right),
\nonumber\\
\mbox{where\ }
 u(t) p(t) := \sum_{i\in I} u_i(t) p_i(t),\notag \\
\label{eq:comp}
1 = \int du(t)  | u(t) \rangle  \langle u(t)|,
\qquad
1 = \int \frac{dp(t)}{(2\pi \hbar)^n} | p(t) \rangle  \langle p(t)|.
\end{gather}
In particular, we have
\begin{gather}
\label{eq:uDp}
\frac{
\langle u(t) |
\hat{D}_i(t)|
p(t) \rangle}
{\langle u(t) |
p(t) \rangle}
=
p_i(t) + w_i(t),
\qquad
w_i(t):=\sum_{j=1}b_{ji}(t) u_j(t).
\end{gather}
Let $\hat O$ be the composition of the operators
in the left hand side of \eqref{eq:QDI7}, namely,
\begin{gather*}
\hat O=
\mathbf{\Phi}_b(\varepsilon_{1}{\hat {\mathsf{D}}}_{k_1}(1))^{\varepsilon_1}
 \rho^*_{k_1,\varepsilon_1}
  \cdots
\mathbf{\Phi}_b(\varepsilon_{L}{\hat {\mathsf{D}}}_{k_L}(L))^{\varepsilon_{L}}
 \rho^*_{k_L,\varepsilon_L}
\nu^*\, (= 1),
\end{gather*}
where $(\varepsilon_1,\dots,\varepsilon_L)$ is the
tropical sign-sequence.
Choose any position eigenvector $| u(1)\rangle$.
Then, set the momentum eigenvector $| \tilde{p}(1) \rangle$ such that
its eigenvalues are given by
\begin{gather}
\label{eq:pw}
\tilde{p}_i(1) = w_i(1):= \sum_{j\in I} b_{ji}(1) u_j(1),
\end{gather}
where the notation  $\tilde{p}(1)$ is used for later
convenience.
The condition \eqref{eq:pw} will be used only at the last stage when we construct
the solution of the saddle point equations in Section
\ref{subsec:sol}.

The main idea of our consideration is to study the semiclassical behavior of the quantum identity by using $q$-$p$ symbols of operators, see for example \cite{Berezin65}.
By Dirac's argument \cite{Dirac58}, the semiclassical limit of a $q$-$p$ symbol of a unitary operator $\mathcal{O}$ is given by the
exponential of the generating function of the canonical transformation, which quantum mechanically corresponds to the unitary inner transformation generated by $\mathcal{O}$. In our case, the $q$-$p$ symbol corresponds to the `$u$-$p$' symbol def\/ined by
\begin{gather*}
F({u}(1),\tilde{p}(1)):=
\frac{\langle {u}(1) | \hat O | \tilde{p} (1) \rangle}
{\langle {u}(1) | \tilde{p} (1) \rangle}.
\end{gather*}
Below we  show that the leading term of $\log F({u}(1),\tilde{p}(1))$
in the limit $\hbar \rightarrow 0$
yields the left hand side of \eqref{eq:DI3}
up to a multiplicative constant.
We know {\em a priori} that its value is 0,
which
yields the right hand side of \eqref{eq:DI3}.

\subsection{Integral expression}
\label{subsec:ie}
By inserting the intermediate complete states \eqref{eq:comp},
we obtain the following integral expression
\begin{gather*}
  F({u}(1),\tilde{p}(1))
=
   (2\pi \hbar)^{-n(2L-1)}
  \int dp(1) d\tilde{p}(2) du(2) dp(2) d\tilde{p}(3) \cdots
d\tilde{p}(L) du(L) dp(L)
\\
\phantom{F({u}(1),\tilde{p}(1))}{}\times
\langle \tilde{p}(1) | {u}(1) \rangle
\frac{\langle {u}(1) |\mathbf{\Phi}_b
(\varepsilon_1 \hat{\mathsf{D}}_{k_1}(1))^{\varepsilon_1}
 |  p(1) \rangle}
{\langle {u}(1) | p(1) \rangle}
\langle {u}(1) | p(1) \rangle
\langle p(1) |  \rho^*_{k_1,\varepsilon_1} | \tilde{p}(2) \rangle
\\
\phantom{F({u}(1),\tilde{p}(1))}{}\times
\langle \tilde{p}(2) | u(2) \rangle
\frac{\langle u(2) |\mathbf{\Phi}_b
(\varepsilon_2 \hat{\mathsf{D}}_{k_2}(2))^{\varepsilon_2}
 |  p(2) \rangle}
{\langle u(2) | p(2) \rangle}
\langle u(2) | p(2) \rangle
\langle p(2) |  \rho^*_{k_2,\varepsilon_2} | \tilde{p}(3) \rangle  \cdots \\
\phantom{F({u}(1),\tilde{p}(1))}{}\times
\langle \tilde{p}(L) | u(L) \rangle
\frac{\langle u(L) |\mathbf{\Phi}_b
(\varepsilon_L \hat{\mathsf{D}}_{k_L}(L))^{\varepsilon_L}
 |  {p}(L)\rangle}
{\langle u(L) |p(L) \rangle}
\langle u(L) | p(L)  \rangle
\langle p(L)
 |  \rho^*_{k_L,\varepsilon_L} \nu^* | \tilde{p}(1) \rangle.
\end{gather*}
The integration over $p(L)$ is done by
\eqref{eq:dtrans}, and it yields the relation
\begin{gather}
\label{eq:sp1}
\tilde{p}_i(1)
=
\begin{cases}
-{p}_{k_L}(L), & \nu(i)= k_L,\\
{p}_{\nu(i)}(L) + [\varepsilon_L b'_{k_L \nu(i)}(L)]_+ {p}_{k_L}(L),
& \nu(i) \neq k_L.
\end{cases}
\end{gather}
Similarly, the integration over $\tilde{p}(t+1)$ ($t=1,\dots,L-1$)
 yields the relation
\begin{gather}
\label{eq:sp2}
\tilde{p}_i(t+1)
=
\begin{cases}
-{p}_{k_t}(t), & i= k_t,\\
{p}_i(t) + [\varepsilon_t b'_{k_ti}(t)]_+ {p}_{k_t}(t),
& i \neq k_t.
\end{cases}
\end{gather}
Thus, $\tilde{p}(t+1)$ is now a dependent
variable of ${p}(t)$ by \eqref{eq:sp2}.

In view of \eqref{eq:wD} it is natural to introduce new dependent
variables
\begin{gather}
\label{eq:ypw}
y_{k_t}(t) =
 \exp (p_{k_t}(t) + w_{k_t}(t)), \qquad
t=1,\dots,L,
\end{gather}
where the notation $y_i(t)$ anticipates the identif\/ication with
classical $y$-variables eventually.
Then, by \eqref{eq:uDp}, we have
\begin{gather}
\label{eq:uDp2}
\frac{
\langle u(t) | \hat{\mathsf{D}}_{k_t} (t) | p(t) \rangle
}
{
\langle u(t)  | p(t) \rangle
}
=\frac{1}{2\pi b} \log y_{k_t}(t),
\end{gather}
and the remaining integration has the following form
\begin{gather*}
 F({u}(1),\tilde{p}(1))
=   (2\pi \hbar)^{-n(L-1)}  \int dp(1)  \cdots
dp(L-1)  du(2)  \cdots du(L)
\\
\phantom{F({u}(1),\tilde{p}(1))=}{}\times
\mathbf{\Phi}_b
\left(
\frac{\log y_{k_1}(1)^{\varepsilon_1}}
{2\pi b}
\right)^{\varepsilon_1}
\exp\left(
\frac{\sqrt{-1}}{\hbar}
{u}(1)(p(1) - \tilde{p}(1))
\right)
\\
\phantom{F({u}(1),\tilde{p}(1))=}{}\times
\mathbf{\Phi}_b
\left(
\frac{\log y_{k_2}(2)^{\varepsilon_2}}
{2\pi b}
\right)^{\varepsilon_2}
\exp\left(
\frac{\sqrt{-1}}{\hbar}
{u}(2)(p(2) - \tilde{p}(2))
\right)
 \cdots \\
\phantom{F({u}(1),\tilde{p}(1))=}{}\times
\mathbf{\Phi}_b
\left(
\frac{\log y_{k_L}(L)^{\varepsilon_L}}
{2\pi b}
\right)^{\varepsilon_L}
\exp\left(
\frac{\sqrt{-1}}{\hbar}
{u}(L)(p(L) - \tilde{p}(L))
\right).
\end{gather*}
Using \eqref{eq:asym2}, we have
\begin{gather}
 F({u}(1),\tilde{p}(1))
 \sim   (2\pi \hbar)^{-n(L-1)}  \int dp(1)   \cdots
dp(L-1)  du(2)  \cdots du(L)
\nonumber\\
\phantom{F({u}(1),\tilde{p}(1)) \sim}{}
\exp\left(
\frac{\sqrt{-1}}{\hbar}
\sum_{t=1}^L
\left\{
\frac{1}{2}
\varepsilon_t
\mathrm{Li}_2(- y_{k_t}(t)^{\varepsilon_t})
+
{u}(t)(p(t) - \tilde{p}(t))
\right\}
\right).\label{eq:int2y}
\end{gather}
To evaluate
 the integral expression \eqref{eq:int2y}
in the semiclassical limit,
we apply the saddle point method.
It consists of three steps.

Step~1. Write the saddle point equations,
that is, the extremum condition of the integrand of~\eqref{eq:int2y}
for the independent variables
$p(1), \dots, p(L-1)$ and $u(2), \dots, u(L)$.

Step~2. Find a  solution of the saddle point equations.

Step~3. Evaluate the integrand at the solution.

\subsection{Saddle point equations}
\label{subsec:spe}

Let us derive the saddle point equations for \eqref{eq:int2y}.
We use the following formulas, which are obtained from
\eqref{eq:L00}, \eqref{eq:wD},
\eqref{eq:ypw}, and~\eqref{eq:uDp2},
\begin{gather*}
\frac{\partial}{\partial p_{i}(t)}
\left(
\frac{1}{2}
\varepsilon_t
\mathrm{Li}_2(- y_{k_t}(t)^{\varepsilon_t})
\right)
 =
\delta_{i k_t}
\log (1+ y_{k_t}(t)^{\varepsilon_t})^{-1/2},\\
\frac{\partial}{\partial u_{i}(t)}
\left(
\frac{1}{2}
\varepsilon_t
\mathrm{Li}_2(- y_{k_t}(t)^{\varepsilon_t})
\right)
 =
- \log (1+ y_{k_t}(t)^{\varepsilon_t})^{-b_{k_ti}(t)/2 }.
\end{gather*}

(a) Extremum conditions  with respect to $u_i(t)$ ($t=2,\dots, L$).

By dif\/ferentiating the integrand of \eqref{eq:int2y} by $u_i(t)$, we have
\begin{gather}
\label{eq:peq}
- \log (1+y_{k_t}(t)^{\varepsilon_t})^{-b_{k_t i}(t)/2} + p_{i}(t) -
\tilde{p}_{i}(t) = 0,
\end{gather}
or, equivalently,
\begin{gather*}
e^{p_{i}(t)}
=
e^{\tilde{p}_{i}(t)}
(1+y_{k_t}(t)^{\varepsilon_t})^{-b_{k_ti}(t)/2}.
\end{gather*}
Combining it with \eqref{eq:sp2}, we also have
\begin{gather}
\label{eq:peq2}
e^{\tilde{p}_{i}(t+1)}
=
\begin{cases}
(e^{\tilde{p}_{k_t}(t)})^{-1},& i = k_t,\\
e^{\tilde{p}_{i}(t)}
(e^{\tilde{p}_{k_t}(t)})^{[\varepsilon_t b_{k_t i} (t)]_+}
(1+y_{k_t}(t)^{\varepsilon_t})^{-b_{k_ti}(t)/2}, &
i \neq k_t.
\end{cases}
\end{gather}

(b) Extremum conditions  with respect to $p_i(t)$ ($t=1,\dots, L-1$).

By dif\/ferentiating the integrand of \eqref{eq:int2y} by $p_i(t)$, we have
\begin{gather}
\label{eq:u2}
 \log (1+y_{k_t}(t)^{\varepsilon_t})^{-1/2} + u_{k_t}(t) -
\sum_{j\in I} [\varepsilon_t b_{k_tj}(t)]_+ u_j(t+1)
+ u_{k_t}(t+1) = 0,
\qquad i = k_t,\!\!\!\!\\
\label{eq:u1}
u_{i}(t) -u_{i}(t+1) = 0, \qquad i\neq k_t,
\end{gather}
or, equivalently,
\begin{gather}
\label{eq:u3}
e^{u_{i}(t+1)}
=
\begin{cases}
(e^{{u}_{k_t}(t)})^{-1}
\prod\limits_{j\in I}
(e^{{u}_{j}(t)})^{[-\varepsilon_t b_{jk_t}(t)]_+}
(1+y_{k_t}(t)^{\varepsilon_t})^{1/2}
& i = k_t,\vspace{1mm}\\
e^{{u}_{i}(t)}, & i\neq k_t.
\end{cases}
\end{gather}
With \eqref{eq:Bmut}, this  also implies
 the following equations for $w_i(u)=
\sum_{j\in I} b_{ji}(t) u_j(t)$
\begin{gather}
\label{eq:weq}
e^{w_{i}(t+1)}
=
\begin{cases}
(e^{w_{k_t}(t)})^{-1},& i = k_t,\\
e^{w_{i}(t)}
(e^{w_{k_t}(t)})^{[\varepsilon_t b_{k_t i} (t)]_+}
(1+y_{k_t}(t)^{\varepsilon_t})^{-b_{k_ti}(t)/2},&
i \neq k_t,
\end{cases}
\end{gather}
which is identical to \eqref{eq:peq2}.

\subsection{Solution}
\label{subsec:sol}
Let us summarize  the relevant variables and their relations schematically
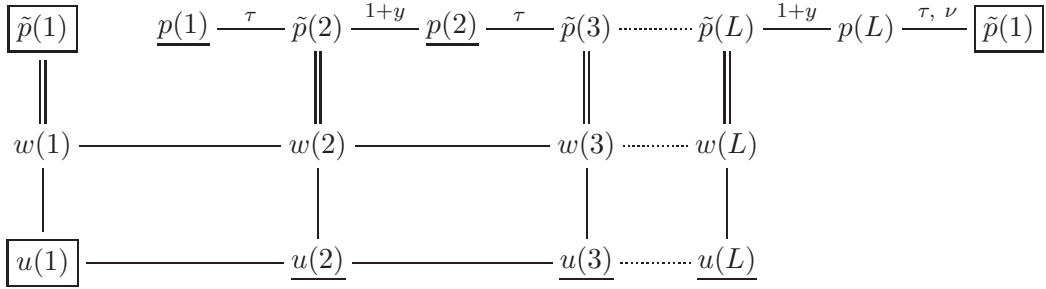
\begin{gather*}
\begin{xy}
\xymatrix{
\framebox{$\tilde{p}(1)$} &
\underline{p(1)}\ar@{-}[r]^{\tau}&\tilde{p}(2)\ar@{-}[r]^{1+y}&
\underline{p(2)}\ar@{-}[r]^{\tau}&\tilde{p}(3)\ar@{.}[r]& \tilde{p}(L)
\ar@{-}[r]^{1+y}&
p(L)\ar@{-}[r]^{\tau,\ \nu}&\framebox{$\tilde{p}(1)$}
\\
w(1)\ar@{-}[rr] \ar@{-}[d]\ar@{=}[u] &&
w(2)\ar@{-}[rr] \ar@{-}[d]\ar@{=}[u]&&
 w(3) \ar@{.}[r] \ar@{-}[d]\ar@{=}[u]
& w(L) \ar@{-}[d]\ar@{=}[u]
\\
\framebox{$u(1)$}\ar@{-}[rr] &&
\underline{u(2)}\ar@{-}[rr] &&
\underline{ u(3)} \ar@{.}[r]
& \underline{u(L)}
}
\end{xy}
\end{gather*}
Here, the framed variables are the initial variables
and the underlined variables are the remaining integration variables
which should be determined to solve the saddle point equations.
This is a~highly complicated systems of equations, but the relevance
to the $y$-seed mutations of \eqref{eq:seedmutseq} is rather clear.
To see it quickly, set
\begin{gather*}
y_i (t):= e^{\tilde{p}_i(t)} e^{w_i(t)}.
\end{gather*}
Note that $\tilde{p}_{i}(t)= {p}_i(t)$ if $i = k_t$ by \eqref{eq:peq},
therefore,
it agrees with the previous def\/inition~\eqref{eq:ypw}.
Then, multiply two identities~\eqref{eq:peq2} and~\eqref{eq:weq}, we have
\begin{gather*}
y_{i}(t+1)
=
\begin{cases}
y_{k_t}(t)^{-1},& i = k_t,\\
y_{i}(t)
y_{k}(t)^{[\varepsilon_t b_{k_t i} (t)]_+}
(1+y_{k_t}(t)^{\varepsilon_t})^{-b_{k_ti}(t)},&
i \neq k_t.
\end{cases}
\end{gather*}
This is nothing but \eqref{eq:coef}.
Furthermore, \eqref{eq:peq2} and \eqref{eq:weq} suggest that
\begin{gather*}
y_i (t)^{1/2}= e^{\tilde{p}_i(t)}= e^{w_i(t)}.
\end{gather*}

Having this observation in mind,
let us describe the construction of the solution more clearly.

$(i)$ ($y$-variables) We have $u_i(1)$ as initial data, from which $w_i(1)$ is
uniquely determined.
Temporarily forgetting~\eqref{eq:ypw},
set $y_i(1)=  e^{2 w_i(1)}$, from which $y_i(t)$ ($t=2,\dots,L$) are
determined by the mutation sequence \eqref{eq:seedmutseq}.

$(ii)$ ($u$-variables)
Set $u_i(t)$ ($t=2,\dots,L$) by \eqref {eq:u1} and \eqref{eq:u2}.
 Then, \eqref{eq:weq} is also satisf\/ied.

$(iii)$ ($p$-variables)
 Set $\tilde{p}_i(t)$ by
$e^{\tilde{p}_i(t)} = y_i (t)^{1/2}$. This forces the relation
$\tilde{p}_i(1)=w_i(1)$, which is guaranteed by the assumption \eqref{eq:pw}.
Then, ${p}_i(t)$ are determined by \eqref{eq:sp1} and \eqref{eq:sp2}.
Since~$\tilde{p}_i(t)$  satisf\/ies \eqref{eq:peq2} by def\/inition,
\eqref{eq:peq} is also satisf\/ied.

$(iv)$ (compatibility) The only thing to be checked is \eqref{eq:ypw}.
Since ${p}_{k_t}(t) = \tilde{p}_{k_t}(t)$ by \eqref{eq:peq},
it is enough to show
\begin{gather}
\label{eq:epw}
e^{\tilde{p}_i(t)} = y_i (t)^{1/2},
\qquad e^{w_i(t)} = y_i (t)^{1/2}.
\end{gather}
The f\/irst equality is by def\/inition.
The second equality is true for $t=1$ by def\/inition.
Then, the rest is shown by~\eqref{eq:coef} and the square of~\eqref{eq:weq}.

Thus, we obtain the desired solution of the saddle point equations.
We do not argue on the uniqueness of the solution here
as stated in Section~\ref{subsec:cq}.

\begin{Remark}
Since \eqref{eq:u3} is the square half of the exchange relation
of the $x$-variables of the corresponding cluster algebras
\cite[Proposition~2.3]{Fock03},
the variable $e^{u_i(t)}$ is regarded as
the square half of the $x$-variable
$x_i(t)$.
\end{Remark}

\subsection{Result}
\label{subsec:res}

 As the f\/inal step, we evaluate
the logarithm of the integrand in~\eqref{eq:int2y}
at the solution of the saddle point equations
in Section~\ref{subsec:sol}.
Using~\eqref{eq:asym2}
and ignoring the common factor, it is given by
\begin{gather}
\label{eq:result1}
\sum_{t=1}^L
\left\{
 \frac{1}{2} \varepsilon_t \mathrm{Li_2} (-y_{k_t}(t)^{\varepsilon_t})
+
\sum_{i\in I}
u_i(t) (p_i(t)-\tilde{p}_i(t))
\right\}.
\end{gather}
Recall that
\begin{gather*}
p_i(t)-\tilde{p}_i(t)
 =\log (1+y_{k_t}(t)^{\varepsilon_t})^{-b_{k_t i}(t)/2},\\
w_i(t) = \frac{1}{2} \log  y_i (t)
\end{gather*}
by \eqref{eq:peq} and \eqref{eq:epw}.
Then, the second term of
\eqref{eq:result1} is rewritten as
\begin{gather}
\sum_{i\in I}
u_i(t) (p_i(t)-\tilde{p}_i(t))
 =
\sum_{i\in I}
 u_i(t)
\log (1+y_{k_t}(t)^{\varepsilon_t})^{-b_{k_t i}(t)/2}\nonumber\\
\hphantom{\sum_{i\in I} u_i(t) (p_i(t)-\tilde{p}_i(t))}{}
=
\frac{1}{2}\left(\sum_{i\in I}
b_{ik_t}(t)  u_i(t)
\right)
\log (1+y_{k_t}(t)^{\varepsilon_t})\nonumber\\
\hphantom{\sum_{i\in I} u_i(t) (p_i(t)-\tilde{p}_i(t))}{}
=
\frac{1}{2}
w_{k_t}(t)
\log (1+y_{k_t}(t)^{\varepsilon_t}) =
\frac{1}{4}
\varepsilon_t
\log y_{k_t}(t)^{\varepsilon_t}
\log (1+y_{k_t}(t)^{\varepsilon_t}).\!\!\!\label{eq:logterm}
\end{gather}
Therefore, by \eqref{eq:LL2}, \eqref{eq:result1} is equal to
\begin{gather*}
- \frac{1}{2}
\sum_{t=1}^L
\varepsilon_t
L\left(
\frac{y_{k_t}(t)^{\varepsilon_t}}
{1+ y_{k_t}(t)^{\varepsilon_t}}
\right),
\end{gather*}
but we know it is 0 from the beginning. This is the classical
dilogarithm identity \eqref{eq:DI3}.

\appendix

\section[Quantum dilogarithm identities in local form for $\mathbf{\Psi}_q(x)$ and their semiclassical limits]{Quantum dilogarithm identities in local form for $\boldsymbol{\mathbf{\Psi}_q(x)}$ \\ and their semiclassical limits}\label{appendixA}

In this section we present the quantum dilogarithm identities
in local form for $\mathbf{\Psi}_q(x)$ with tropical sign-sequence.
Then, we derive the classical dilogarithm identities from them
in the semiclassical
limits. The treatment is parallel to the one in
Sections~\ref{sec:phi} and~\ref{sec:from} with slight complication.

\subsection[Representation of quantum $y$-variables]{Representation of quantum $\boldsymbol{y}$-variables}

We consider a representation of
quantum $y$-variables as dif\/ferential operators
which are quite similar to the one in Section~\ref{subsec:qy} but slightly
dif\/ferent.

Throughout the section, let $\hbar$ be a positive real number,
and $\lambda$ be a complex number such that
\begin{gather*}
\mathrm{Im}\, \lambda^2 > 0.
\end{gather*}
We reset
\begin{gather}
\label{eq:qa}
q=e^{\lambda^2\hbar \sqrt{-1}}.
\end{gather}
By the assumption, we have $|q|<1$.
Compare with $q$ in \eqref{eq:hbar},
where $|q|=1$ when $b$ is real.
This dif\/ference is due to the
fact that
$\mathbf{\Psi}_q(x)$ is convergent only for $|q|<1$,
while   $\mathbf{\Phi}_b(z)$ is well-def\/ined also for $|q|=1$.
The phase $\lambda$
is the main dif\/ference between the two cases
and the source of extra complication for $\mathbf{\Psi}_q(x)$
which persists throughout the section.

The asymptotic property \eqref{eq:asym4} is written as
\begin{gather}
\label{eq:asym2a}
\mathbf{\Psi}_q(x) \sim
\exp \left( \frac{\sqrt{-1}}{\lambda^2\hbar}\frac{1}{2}\mathrm{Li}_2(-x)
\right),
\qquad \hbar \rightarrow 0.
\end{gather}
Because of $\lambda$, the argument $x$ of the dilogarithms
$\mathrm{Li}_2(x)$ and $L(x)$ eventually take values in $\mathbb{C}$
 in the semiclassical limit.
They are def\/ined by analytic continuation of
\eqref{eq:L00} and \eqref{eq:L0} along the integration path.
To avoid the ambiguity of the branches,
we assume that {\em $\mathrm{Im}\, \lambda$ is sufficiently small}
(or, $q$ is suf\/f\/iciently close to the unit circle $|q|=1$) so that
the resulting argument $x$ in this section is in a neighborhood of  the interval
$(-\infty, 1]$ for $\mathrm{Li}_2(x)$ or
$[0,1]$ for $L(x)$.

To any quantum  $y$-seed $(B',Y')$ of  $(B,Y)$
 we associate operators $\hat u'=(\hat u'_i)_{i\in I}$
 and $\hat p'=(\hat p'_i)_{i\in I}$,
and the local coordinates $\{\langle u'\vert\}_{u'\in\mathbb{R}^I}$
 as in \eqref{eq:ud2} and \eqref{eq:lc1}.

We reset ${\mathsf{\hat{Y}}}'_i$ in \eqref{eq:wD} as
\begin{gather}
\label{eq:wDa}
\hat w'_i =
\sum_{j\in I} b'_{ji} \hat u'_j,
\qquad
\hat D'_i =
 \hat p'_i
 + \hat w'_i,
\qquad
{\mathsf{\hat{Y}}}'_i= \exp(\lambda \hat D'_i).
\end{gather}
The relations in~\eqref{eq:dd1} still hold
with $q$ in~\eqref{eq:qa},
and we have a representation of
$\mathbb{T}(B,\mathsf{Y})$ of~\eqref{eq:Yrep1}.

Let $(B',Y')$ and $(B'',Y'')$ be a pair of quantum
$y$-seeds of $(B,Y)$ such that
$(B'',Y'')=\mu_k (B',Y')$.
Let $\rho_{k,\varepsilon}$ be the map
in \eqref{eq:utrans}.
Then, repeating the argument in Section \ref{subsec:mut2},
we obtain
\begin{gather*}
\mathsf{Ad}(
\rho^*_{k,\varepsilon}
)
(\mathsf{\hat{Y}}''_i) =
\begin{cases}
\mathsf{\hat{Y}}'_k{}^{-1}, & i= k,\\
\mathsf{\hat{Y}}'^{e_i +  [\varepsilon b'_{ki}]_+ e_k},
& i \neq k.
\end{cases}
\end{gather*}

\subsection[Quantum dilogarithm identities in local form for $\mathbf{\Psi}_q(x)$]{Quantum dilogarithm identities in local form for $\boldsymbol{\mathbf{\Psi}_q(x)}$}

Under the same assumption and notation for Theorem~\ref{thm:QDI5},
we obtain the counterpart of
Theorem~\ref{thm:QDI5} for~$\mathbf{\Psi}_q(x)$
by repeating its proof.

\begin{Theorem}
\label{thm:QDI5a}
For the tropical
sign-sequence $
{\vec{\varepsilon}}=
(\varepsilon_{1}, \dots, \varepsilon_{L})$,
the following identity holds
\begin{align}
\label{eq:QDI7a}
\mathbf{\Psi}_q(\mathsf{\hat{Y}}_{k_1}(1)^{\varepsilon_{1}})^{\varepsilon_1}
 \rho^*_{k_1,\varepsilon_1}
  \cdots
\mathbf{\Psi}_q(\mathsf{\hat{Y}}_{k_L}(L)^{\varepsilon_{L}})^{\varepsilon_{L}}
 \rho^*_{k_L,\varepsilon_L}
\nu^*
=1.
\end{align}
\end{Theorem}

Let $\hat O$ be the composition of the operators
in the left hand side of \eqref{eq:QDI7a}.
Again, choose any position eigenvector $| u(1)\rangle$
and set the momentum eigenvector $| \tilde{p}(1) \rangle$ by
\eqref{eq:pw}.
Set
\begin{gather*}
F({u}(1),\tilde{p}(1)):=
\frac{\langle {u}(1) | \hat O | \tilde{p} (1) \rangle}
{\langle {u}(1) | \tilde{p} (1) \rangle}.
\end{gather*}
Below we  show that the leading term of $\log F({u}(1),\tilde{p}(1))$
in the limit $\hbar \rightarrow 0$
yields the left hand side of~\eqref{eq:DI3}
up to a multiplicative constant.

\subsection{Integral expression}

Repeating the argument in Section~\ref{subsec:ie},
we obtain the following integral expression
\begin{gather*}
 F({u}(1),\tilde{p}(1))
=  (2\pi \hbar)^{-n(L-1)}  \int dp(1)   \cdots
dp(L-1)  du(2)  \cdots du(L)
\\
\phantom{F({u}(1),\tilde{p}(1))=}{}
\times\mathbf{\Psi}_q
\left(
{ y_{k_1}(1)^{\varepsilon_1}}
\right)^{\varepsilon_1}
\exp\left(
\frac{\sqrt{-1}}{\hbar}
{u}(1)(p(1) - \tilde{p}(1))
\right)
\\
\phantom{F({u}(1),\tilde{p}(1))=}{}
\times\mathbf{\Psi}_q
\left(
{ y_{k_2}(2)^{\varepsilon_2}}
\right)^{\varepsilon_2}
\exp\left(
\frac{\sqrt{-1}}{\hbar}
{u}(2)(p(2) - \tilde{p}(2))
\right)
 \cdots \\
\phantom{F({u}(1),\tilde{p}(1))=}{}
\times\mathbf{\Psi}_q
\left(
{ y_{k_L}(L)^{\varepsilon_L}}
\right)^{\varepsilon_L}
\exp\left(
\frac{\sqrt{-1}}{\hbar}
{u}(L)(p(L) - \tilde{p}(L))
\right),
 \end{gather*}
where
 $\tilde{p}(t)$ is the one in Section \ref{subsec:ie},
while we reset
\begin{gather}
\label{eq:ypwa}
y_{k_t}(t) =
 \exp \left(\lambda (p_{k_t}(t) + w_{k_t}(t))\right).
\end{gather}
Using \eqref{eq:asym2a}, we have
\begin{gather}
 F({u}(1),\tilde{p}(1))
\sim   (2\pi \hbar)^{-n(L-1)}  \int dp(1)   \cdots
dp(L-1)  du(2)  \cdots du(L)
\nonumber\\
\phantom{F({u}(1),\tilde{p}(1))\sim}{}
\times\exp\left(
\frac{\sqrt{-1}}{\hbar}
\sum_{t=1}^L
\left\{
\frac{1}{2\lambda^2}
\varepsilon_t
\mathrm{Li}_2(- y_{k_t}(t)^{\varepsilon_t})
+
{u}(t)(p(t) - \tilde{p}(t))
\right\}
\right).\label{eq:int2ay}
\end{gather}

\subsection{Saddle point equations}

The saddle point equations for \eqref{eq:int2ay}
are obtained in the same manner as in Section~\ref{subsec:spe}.
We use the following formulas, which are obtained from
\eqref{eq:L00},   \eqref{eq:wDa},  and
\eqref{eq:ypwa}
\begin{gather*}
\frac{\partial}{\partial p_{i}(t)}
\left(
\frac{1}{2\lambda^2}
\varepsilon_t
\mathrm{Li}_2(- y_{k_t}(t)^{\varepsilon_t})
\right)
 =
\delta_{i k_t}
\frac{1}{\lambda}
\log (1+ y_{k_t}(t)^{\varepsilon_t})^{-1/2},\\
\frac{\partial}{\partial u_{i}(t)}
\left(
\frac{1}{2\lambda^2}
\varepsilon_t
\mathrm{Li}_2(- y_{k_t}(t)^{\varepsilon_t})
\right)
 =
- \frac{1}{\lambda}
 \log (1+ y_{k_t}(t)^{\varepsilon_t})^{-b_{k_ti}(t)/2 }.
\end{gather*}

(a) Extremum conditions  with respect to $u_i(t)$ ($t=2,\dots, L$).

By dif\/ferentiating the integrand of \eqref{eq:int2ay} by $u_i(t)$, we have
\begin{gather}
\label{eq:peqa}
- \frac{1}{\lambda}
\log (1+y_{k_t}(t)^{\varepsilon_t})^{-b_{k_t i}(t)/2} + p_{i}(t) -
\tilde{p}_{i}(t) = 0.
\end{gather}
Combining it with \eqref{eq:sp2}, we also have
\begin{gather}
\label{eq:peq2a}
 e^{\lambda\tilde{p}_{i}(t+1)}
=
\begin{cases}
(e^{\lambda\tilde{p}_{k_t}(t)})^{-1},& i = k_t,\\
 e^{\lambda\tilde{p}_{i}(t)}
(e^{\lambda
\tilde{p}_{k_t}(t)})^{[\varepsilon_t b_{k_t i} (t)]_+}
(1+y_{k_t}(t)^{\varepsilon_t})^{-b_{k_ti}(t)/2},&
i \neq k_t.
\end{cases}
\end{gather}

(b) Extremum conditions  with respect to $p_i(t)$ ($t=1,\dots, L-1$).

By dif\/ferentiating the integrand of \eqref{eq:int2ay} by $p_i(t)$, we have
\begin{gather}
\frac{1}{\lambda} \log (1+y_{k_t}(t)^{\varepsilon_t})^{-1/2} + u_{k_t}(t) \nonumber\\
\qquad{} -
\sum_{j\in I} [\varepsilon_t b_{k_tj}(t)]_+ u_j(t+1)
+ u_{k_t}(t+1) = 0,
\quad i = k_t,\label{eq:u2a}\\
\label{eq:u1a}
u_{i}(t) -u_{i}(t+1) = 0, \qquad i\neq k_t.
\end{gather}
With \eqref{eq:Bmut}, this  also implies
 the following equations for $w_i(u)=
\sum_{j\in I} b_{ji}(t) u_j(t)$.
\begin{gather}
\label{eq:weqa}
 e^{\lambda w_{i}(t+1)}
=
\begin{cases}
( e^{\lambda w_{k_t}(t)})^{-1},& i = k_t,\\
 e^{\lambda w_{i}(t)}
( e^{\lambda w_{k_t}(t)})^{[\varepsilon_t b_{k_t i} (t)]_+}
(1+y_{k_t}(t)^{\varepsilon_t})^{-b_{k_ti}(t)/2},&
i \neq k_t.
\end{cases}
\end{gather}

\subsection{Solution}
\label{subsec:sola}

The (complex) solution of the saddle point equations is
constructed in the same manner as in
Section \ref{subsec:sol}
and  given as follows.

$(i)$ ($y$-variables) We have $u_i(1)$ as initial data, from which $w_i(1)$ is
uniquely determined.
Temporarily forgetting \eqref{eq:ypwa},
set $y_i(1)=   e^{2 \lambda
 w_i(1)}$, from which $y_i(t)$ ($t=2,\dots,L$) are
determined by the mutation sequence \eqref{eq:seedmutseq}.

$(ii)$ ($u$-variables)
Set $u_i(t)$ ($t=2,\dots,L$) by \eqref {eq:u1a} and~\eqref{eq:u2a}.
 Then, \eqref{eq:weqa} is also satisf\/ied.

$(iii)$ ($p$-variables)
 Set $\tilde{p}_i(t)$ by $e^{\lambda
\tilde{p}_i(t)} =y_i (t)^{1/2}$. This forces the relation
$\tilde{p}_i(1)=w_i(1)$, which is guaranteed by the assumption \eqref{eq:pw}.
Then, ${p}_i(t)$ are determined by \eqref{eq:sp1} and \eqref{eq:sp2}.
Since~$\tilde{p}_i(t)$  satisf\/ies \eqref{eq:peq2a} by def\/inition,
\eqref{eq:peqa} is also satisf\/ied.

\subsection{Result}

 The evaluation of
the logarithm of the integrand in \eqref{eq:int2ay}
at the solution of the saddle point equations
is done in the same manner as in
Section~\ref{subsec:res}.
Using~\eqref{eq:asym2a}
and ignoring the common factor, it is given by
\begin{gather*}
\sum_{t=1}^L
\left\{
 \frac{1}{2} \varepsilon_t \mathrm{Li_2} (-y_{k_t}(t)^{\varepsilon_t})
+
\lambda^2
\sum_{i\in I}
u_i(t) (p_i(t)-\tilde{p}_i(t))
\right\}.
\end{gather*}
Recall that
\begin{gather*}
\lambda(p_i(t)-\tilde{p}_i(t))
 =\log (1+y_{k_t}(t)^{\varepsilon_t})^{-b_{k_t i}(t)/2},\\
\lambda w_i(t) = \frac{1}{2} \log  y_i (t).
\end{gather*}
Then, repeating the calculation in \eqref{eq:logterm},
 we  obtain the classical
dilogarithm identity \eqref{eq:DI3}
with complex argument.
Taking the  limit $\lambda \rightarrow  1$ further,
we recover the identity \eqref{eq:DI3} with real argument.

\subsection*{Acknowledgments}
We thank
Vladimir V.~Bazhanov, Ludwig D.~Faddeev, Kentaro Nagao,
Boris Pioline, and Andrei Zelevinsky for very useful discussions and comments.
We especially thank Alexander Yu.~Volkov for making his result
in~\cite{Volkov11b} available to us prior to the publication.

\pdfbookmark[1]{References}{ref}
\LastPageEnding

\end{document}